\documentclass [12pt] {amsart}
\usepackage{amsmath, amsfonts, latexsym, array, amssymb, }
\usepackage{psfrag}
\usepackage[all]{xypic}
\usepackage{graphicx}
\newtheorem{thm}{Theorem}[section]
\newtheorem{cor}[thm]{Corollary}
\newtheorem{lem}[thm]{Lemma}
\newtheorem{prop}[thm]{Proposition}
\newtheorem{cla}[thm]{Claim}
\newtheorem{defn}[thm]{Definition}
\newtheorem{rem}[thm]{Remark}

\begin{document}

\title[Resolving extensions] {Resolving extensions of finitely presented systems}
\author{Todd Fisher}
\address{Department of Mathematics Brigham Young University, Provo, UT 84602}
\email{tfisher@math.byu.edu}
\thanks{}

\date{November 2, 2007}
\commby{}

\begin{abstract}
In this paper we extend certain central results of zero dimensional systems to higher dimensions.  The first main result shows that if $(Y,f)$ is a finitely presented system,
then there exists a Smale space $(X,F)$ and a $u$-resolving factor map
$\pi_+: X\rightarrow Y$.   If the finitely presented system is transitive, then we show there is a canonical minimal $u$-resolving Smale space extension.  
Additionally, we show that any finite-to-one factor map between transitive finitely presented systems lifts
through $u$-resolving maps to an $s$-resolving map.
 \end{abstract}

\maketitle

\section{Introduction}

One cornerstone of the study of dynamical systems is the theory of hyperbolic dynamics introduced by  Smale and Anosov in the 1960s.  For compact spaces the property of hyperbolicity, in general, produces highly nontrivial and interesting dynamics.
The best understood hyperbolic sets are those that are locally maximal (or isolated).  For many years it was asked if every hyperbolic set can be contained in a locally maximal hyperbolic set: in~\cite{Fis1}  it is shown that this is not the case for  any manifold with dimension greater than one.  

This paper is in part an investigation into hyperbolic sets that are not locally maximal.  If a hyperbolic set is not locally maximal one of the next properties one considers is the existence of a Markov partition, which is  roughly a decomposition of the hyperbolic set into dynamically defined rectangles.
In~\cite{Fis1} it is shown that any hyperbolic set can be extended to one with a Markov partition.

The essential topological
structure of a locally maximal hyperbolic set is captured
with the notion of a Smale space, introduced by   Ruelle~\cite{Rue}
and simplified by Fried ~\cite{Fri}.  A {\bf Smale space} is an expansive system with canonical coordinates (or a local product structure).

Fried defined  a {\bf finitely presented}  dynamical system as
an expansive system which is a factor  of a
shift of finite type. The finitely presented systems contain
the Smale spaces and share a great deal of their structure; in
particular, they are precisely the expansive systems that admit Markov
partitions~\cite{Fri}.

In this paper we extend certain central results of zero dimensional finitely presented systems (sofic shifts)  
to higher dimensional finitely presented systems. This investigation
is motivated in part by hyperbolic sets that need not be
locally maximal, and also as a generalization to
higher dimension of the symbolic viewpoint. More specifically, the present work looks at resolving maps from Smale spaces to finitely presented
systems; these are absolutely central to the zero-dimensional
theory. 
A factor map from a
space $(X, f )$ to $(Y, g)$ is {\bf $u$-resolving} ({\bf $s$-resolving}) if it is
injective on unstable (stable) sets.

Every sofic shift has a canonical extension to a closely related
 SFT such that the
factor map defining the extension is $u$-resolving.
If the sofic shift is transitive, 
then the extension can be chosen to be transitive.  Furthermore, in this setting there is a  
 minimal extension in the class of 
$u$-resolving extensions of the sofic shift to transitive 
SFTs (i.e., every such extension factors through the 
minimal one). 
In dimension zero these covers permit the study of sofic shifts by the considerably more accessible 
SFTs. 
It is natural to seek such 
covers of finitely presented systems by Smale spaces in 
positive dimension.

\begin{thm}\label{thm1}
If $Y$ is a compact metric space and $(Y,f)$ is finitely presented,
then there exists a Smale space $(X,F)$ and a $u$-resolving factor map
$\pi^+: X\rightarrow Y$.  Furthermore, if $(Y,f)$ is transitive, then 
there exists a transitive Smale space $(X,F)$ 
and a $u$-resolving one-to-one almost everywhere  map $\pi_+:X\rightarrow Y$ (a residual set of points exists with a unique preimage).
\end{thm}

Our proofs use a (highly nontrivial) symbolic coding to
construct the Smale space extensions.
Fried gives a different proof of Theorem~\ref{thm1} in~\cite{Fri}.
One aspect of this construction is that we are 
able to avoid a difficult part of Fried's argument 
[page 494, lines 28-29]. 

In general, there is no minimal extension to a Smale space for a
finitely presented system, see~\cite{Wil1, Wil2}.  However, if $(Y,f)$ is finitely
presented and
transitive we have the following result.

\begin{thm}\label{thm1a}
If $(Y,f)$ is a transitive finitely presented system, then there
exists a minimal extension among all transitive Smale spaces with
u-resolving factor maps onto $Y$.
\end{thm}

Another important class of zero dimensional systems is 
almost of finite type (AFT). 
A subshift is {\bf AFT} if it is the image of a transitive SFT 
via a factor map which is one-to-one on a nonempty open set.  
The AFT subshifts form  the one naturally
distinguished good class of sofic shifts. In particular, 
Boyle, 
Kitchens, and Marcus~\cite{BKM} show that a transitive 
AFT shift,  has a canonical minimal extension to a transitive SFT 
where the extension is minimal in the class of {\it all} extensions 
to transitive SFTs. This extension is $u$-resolving, 
$s$-resolving and one-to-one almost everywhere.

We 
extend this result to positive dimensions.
An expansive system is an {\bf almost Smale space} if it is the image of a transitive
Smale space 
via a factor map which is 1-1 on a nonempty open set.

\begin{thm}\label{thm2}
Let $(Y,f)$ be an almost Smale system, $(X,g)$ be a transitive Smale space,
and
$\theta:X\rightarrow Y$ be 1-1 almost everywhere, $u$-resolving, and
$s$-resolving.  Let $(X', g')$ be an irreducible Smale space and
$\phi:X'\rightarrow Y$ a factor map.  Then $\phi$ factors through $\theta$.
\end{thm}

Approaching this subject with an eye to certain  
$C^*$-algebras 
 Putnam~\cite{Put} proved a surprising result: many factor maps lift
through $u$-resolving maps to $s$-resolving maps.
The next result extends
Putnam's theorem to the case of finitely presented systems.

\begin{thm} \label{thm2a}
Let $(X,f)$ and $(Y,g)$ be transitive and finitely presented, and
$\pi: X\rightarrow Y$ a finite-to-one factor map. Then there exist
Smale spaces $(\tilde{X},\tilde{f})$ and $(\tilde{Y},\tilde{g})$ and factor
maps $\phi$, $\theta$, and $\tilde{\pi}$, such that $\phi$ and
$\theta$ are $u$-resolving, $\tilde{\pi}$ is $s$-resolving and the
following diagram commutes:

$$
\xymatrix{
\tilde{X} \ar[d]_{\phi}\ar[r]^{\tilde{\pi}} & \tilde{Y}\ar[d]^{\theta}\\
X\ar[r]_{\pi} &  Y} $$
\end{thm}

We remark that in the zero dimensional case,
Boyle~\cite{Boy2} obtained the above result for sofic
systems,
and even provided a meaningful generalization when the map
is infinite-to-one.  Boyle also establishes a canonical
mapping property in dimension zero.  From Theorem~\ref{thm2a} we have the following corollary.

\begin{cor}
If $(Y,f)$ is transitive and finitely presented, then there exists an SFT
$(\Sigma, \sigma)$, a Smale space $(X,F)$, an s-resolving factor
map $\alpha: \Sigma\rightarrow X$, and a u-resolving factor map
$\beta: X\rightarrow Y$.
\end{cor}

The paper proceeds as follows: In Section 2 we provide some background results and definitions.  In Section 3 we give preliminary results that will be used throughout the rest of the paper.  These results give some basic facts on factor maps and resolving maps.  In Section 4 we provide a proof of Theorem~\ref{thm1}.  In Section 5 we give a proof of Theorem~\ref{thm1a}.  In Section 6 we give a proof of Theorem~\ref{thm2}.  Lastly, in Section 7 we provide a proof of Theorem~\ref{thm2a} and give a meaningful extension of a magic word to the case of a finite-to-one factor map from an SFT to a Smale space.\\

\noindent{\bf Acknowledgement:}  The author would like to thank Mike Boyle for pointing out these problems and offering helpful advice and remarks.

\section{Background}

Before proceeding we review some
useful concepts.  
Throughout the presented work we let $X$ be a compact metric space and $f$ be a homeomorphism of $X$.  We now review expansive dynamics.

\begin{defn} Let $X$ be a compact topological space with metric $d(\cdot,\cdot)$.
A homeomorphism $f:X\rightarrow X$ is {\bf expansive} if there
exists a constant $c>0$ such that for all $p,q\in X$, where $p\neq
q$, there is an $n\in\mathbb{Z}$ with $d(f^n(p),f^n(q))>c$.
\end{defn}

The constant $c>0$ is called an {\bf expansive constant}.  The
definition does not depend on the choice of metric, although the
constant $c$ may.

For $\epsilon>0$ and $x\in X$ the \textbf{ $\epsilon$-stable set} is
$$
W_{\epsilon}^{s}(x)=\{ y\in X\, |\textrm{ for all } n\geq 0,
d(f^{n}(x), f^{n}(y))<\epsilon\},$$ and the
\textbf{$\epsilon$-unstable set} is
$$W_{\epsilon}^{u}(x)=\{ y\in
X\, |\textrm{ for all } n\geq 0, d(f^{-n}(x),
f^{-n}(y))<\epsilon\}.$$

For $x\in X$ and $f:X\rightarrow X$ an expansive
homeomorphism the {\bf
stable set} is
$$W^s(x)=\{y\in X\, |\,
\lim_{n\rightarrow\infty}d(f^{n}(x),f^{n}(y))=0\}$$ and the {\bf unstable set} is
$$W^u(x)=\{y\in X\, |\,
\lim_{n\rightarrow\infty}d(f^{-n}(x),f^{-n}(y))=0\}.$$ 

In~\cite{Fri} it is shown that for any expansive system $(X,f)$
there exists an {\bf adapted metric} $d(\cdot,\cdot)$ and
constants $\epsilon>0$ and $\lambda\in(0,1)$ such that
$$d(f(x),f(y))\leq\lambda d(x,y)\textrm{ for all }y\in
W^s_{\epsilon}(x)\textrm{ and}$$
$$d(f^{-1}(x),f^{-1}(z))\leq\lambda d(x,z)\textrm{ for all }z\in
W^u_{\epsilon}(x).$$
The following is a standard result for expansive systems. 

\begin{lem}
If
$f:X\rightarrow X$ is expansive with expansive constant $c$  and
$\epsilon<c/2$, then for any $x,y\in X$ the intersection
$W^s_{\epsilon}(x)\cap W^u_{\epsilon}(y)$ consists of at most one
point.
\end{lem}

A system $(X,f)$ has {\bf canonical coordinates} if for all sufficiently small $\epsilon>0$ there exists a $ \delta>0$ such that if $x,y\in X$ and $d(x,y)<\delta$, then $W^s_{\epsilon}(x)\cap W^u_{\epsilon}(y)$ consists of one point.

A sequence $\{x_n\}_{j_1}^{j_2}$, where $-\infty\leq j_1<j_2\leq\infty$, is an
{\bf $\epsilon$-pseudo orbit} if $d(f(x_{j}),x_{j+1})<\epsilon$ for all $j_1\leq j< j_2$ if $j_1>-\infty$ or for all $j<j_2$ if $j_1=-\infty$.   For a given $\delta,\epsilon>0$ an $\epsilon$-pseudo orbit is said to be {\bf $\delta$-shadowed} by $x\in X$ if $d(f^i(x),x_i)<\delta$ for all $j_1\leq i\leq j_2$.  A map $f$ has the {\bf pseudo orbit tracing property} (POTP) if for every $\delta>0$ there is $\epsilon>0$ such that every $\epsilon$-pseudo orbit is $\delta$-shadowed by a point $x\in X$.    
In~\cite{Omb} it is shown that a space $(X,f)$ is a Smale space if and only if it has the 
pseudo-orbit tracing property.  

A compact metric space $X$ is {\bf topologically transitive} for a map $f$ if for any open sets $U$ and $V$ in $X$ there is $n\in\mathbb{N}$ such that $f^{-n}(U)\cap V\neq\emptyset$.  A space $X$ is {\bf topologically mixing} for $f$ if for any nonempty open sets $U$ and $V$ in $X$ there is $N\in\mathbb{N}$ such that $f^{-n}(U)\cap V\neq\emptyset$ for all $n\geq N$. 
A point $x\in X$ is a {\bf chain recurrent} point if for all $\epsilon>0$ there is an $\epsilon$-pseudo orbit from $x$ to $x$.  The {\bf chain recurrent set} denoted 
$R(f)$ consists of all chain recurrent points.

The following theorem is a useful description of the chain recurrent set of a Smale space.  The theorem can be found in~\cite[p. 101-102]{AH}, however, the terminology is quite different.

\begin{thm}\label{spectral}(Spectral Decomposition Theorem) Let $(X,f)$ be a Smale space.  Then there exists
a finite collection of disjoint compact invariant sets $B_i$ for $1\leq i\leq l$ such that 
\begin{itemize}
\item $f|_{B_i}:B_i\rightarrow B_i$ is transitive for all $1\leq i\leq l$, 
\item $R(f)=\bigcup_{1\leq i\leq l}B_i$, and
\item $(B_i, f|_{B_i})$ is a Smale space.
\end{itemize}
\end{thm}

We now review some properties of shift spaces.
Let $\mathcal{A}$ be a finite set of elements and $\mathcal{A}^{\mathbb{Z}}$ be the collection of all bi-infinite sequences of symbols from $\mathcal{A}$.  A {\bf block} (or {\bf word}) over $\mathcal{A}$ is a finite sequence of symbols in $\mathcal{A}$.  If $x\in\mathcal{A}^{\mathbb{Z}}$ and $w$ is a block over $\mathcal{A}$, we say $w$ occurs in $x$ if there exist $i,j\in\mathbb{Z}$ such that $w=x[i,j]$.  Let $\mathcal{F}$ be a collection of blocks over $\mathcal{A}$, called the {\bf forbidden blocks}.  Define $X_{\mathcal{F}}$ to be the subset of $\mathcal{A}^{\mathbb{Z}}$ which does not contain any block in $\mathcal{F}$.

A {\bf shift space} $(\Sigma, \sigma)$
consists of a set $\Sigma\subset\mathcal{A}^{\mathbb{Z}}$ such that $\Sigma=X_{\mathcal{F}}$ for some collection $\mathcal{F}$ of forbidden blocks in $\mathcal{A}$ and
$\sigma$ is the shift map, where $(\sigma s)_i=s_{i+1}$ for any
$s\in\Sigma$. 

\begin{defn}
A {\bf subshift of finite type} (SFT) is a shift space $X$ having the form $X_{\mathcal{F}}$ for some finite set $\mathcal{F}$ of forbidden blocks. 
\end{defn} 

For a point $t\in\Sigma$ the
{\bf local unstable set of $s$} is denoted
$$W^u_{\mathrm{loc}}(t)=\{t'\in\Sigma\, |\, t_j=t'_j\textrm{ for
all }j\leq 0\}$$ and the {\bf local stable set of $s$} is denoted
$$W^s_{\mathrm{loc}}(t)=\{t'\in\Sigma\, |\, t_j=t'_j\textrm{ for
all }j\geq 0\}.$$
For $t,t'\in\Sigma$ such that $t_0=t'_0$ define the map
$$[t,t']=\left\{
\begin{array}{ll}
t_i & \textrm{for }i\geq0\\
t'_i & \textrm{for }i\leq 0
\end{array}\right. .$$
So $[t,t']= W^s_{\mathrm{loc}}(t)\cap W^u_{\mathrm{loc}}(t')$.

A {\bf factor map} from $(X,F)$ to $(Y,f)$ is a surjective semi-conjugacy.
For $(Y,f)$ a factor of an SFT $(\Sigma, \sigma)$ by a factor map $\pi:\Sigma\rightarrow Y$ we define a relation
$E_{\pi}\subset\Sigma\times\Sigma$ by $(s,t)\in E_{\pi} $ if
$\pi(s)=\pi(t)$. It is clear that $E_{\pi}$ is an equivalence
relation.  The following Lemma provides a useful characterization of
$E_{\pi}$.

\begin{lem}\label{quotientlem}~\cite{Fri} Let $(\Sigma, \sigma)$ be an SFT and $(Y,f)$ a
factor of $(\Sigma, \sigma)$ such that $Y\simeq\Sigma/ E_{\pi}$ for
an equivalence relation $E_{\pi}$.  Then $(Y,f)$ is expansive if and only
if $E_{\pi}$ is an SFT.
\end{lem}

Let $(Y,f)$ be expansive and $\epsilon<c/2$ for $c$ an expansive constant of $(Y,f)$.
Following Fried in~\cite{Fri} we define
$$D_{\epsilon}=\{(x,y)\in Y\times Y\, |\, W^s_{\epsilon}(x)\textrm{
meets }W^u_{\epsilon}(y)\}$$ and $[\cdot,
\cdot]:D_{\epsilon}\rightarrow Y$ so that
$[x,y]=W^s_{\epsilon}(x)\cap W^u_{\epsilon}(y)$. It follows that
$[\cdot, \cdot]$ is continuous.

\begin{defn} A {\bf rectangle} is a closed set $R\subset Y$ such that $R\times R\subset D_{\epsilon}$.
\end{defn}

For $R$ a rectangle and $x\in R$ denote
$$W^s(x,R)=R\cap
W^s_{\epsilon}(x)\textrm{ and }W^u(x,R)=R\cap W^u_{\epsilon}(x).$$
Let $h_x:R\rightarrow W^s_{\epsilon}(x,R)\times W^u_{\epsilon}(x,R)$
so that $h_x(y)=([x,y], [y,x])$.  One easily checks that this is a
homeomorphism with $h_x^{-1}(y,z)=[y,z]$.
A rectangle $R$ is \textbf{proper} if
$R=\overline{\mathrm{int}R}$. 

\begin{defn} Let $(Y,f)$ be expansive with constant $c>0$ and $0<\epsilon<c/2$.  A finite cover $\mathcal{R}$ of $Y$ by
proper rectangles with diameter$(R)<\epsilon$ for any
$R\in\mathcal{R}$ is a {\bf Markov Partition} if 
$R_i,R_j\in\mathcal{R}$,  $x\in\mathrm{int}R_i$, and
$f(x)\in\mathrm{int}R_j$, then
\begin{itemize}
\item
$f(W^s(x,R_i))\subset R_j$ and 
$f^{-1}(W^u(f(x),R_j))\subset R_i$, and
\item $\mathrm{int}(R_i)\cap\mathrm{int}(R_j)=\emptyset$ if $i\neq j$.
\end{itemize}
\end{defn}

Let $(Y,f)$ be a finitely presented system with expansive
constant $c$ and assume there exists an SFT $(\Sigma, \sigma)$ with a factor map $\pi:
\Sigma\rightarrow Y$.  Fix $\epsilon<c/2$.  After passing to a higher block presentation we may assume that for each
$j\in\mathcal{A}(\Sigma)$ the cylinder set $C_j=\{s\in\Sigma\, |
s_0=j\}$ has an image $R_j=\pi(C_j)$ that is a rectangle in $Y$.
Furthermore, if $i, j\in\mathcal{A}(\Sigma)$ with an allowed
transition from $i$ to $j$, $y\in R_i$, and $f(y)\in R_j$, then
$$f(W^s(y, R_i))\subset W^s(f(y), R_j)\textrm{ and
}f^{-1}(W^u(f(y),R_j)\subset W^u(y, R_i).$$
This is called the {\bf Markov property}.

For $R_i$ and $R_j$ rectangles the {\bf unstable $j$ boundary of $R_i$} is
$$\partial_u^{j}R_i=\{x\in R_i\,|\,x=\lim x_i, x_i\in W^s_{\epsilon}(x)\cap R_j- R_i\}.$$  Similarly, we can define the {\bf stable $j$ boundary of $R_i$}
as
$$\partial_s^{j}R_i=\{x\in R_i\,|\,x=\lim x_i, x_i\in W^u_{\epsilon}(x)\cap R_j- R_i\}.$$
For a finitely presented system the {\bf boundary of a rectangle} is 
$$\partial R_i=\bigcup_{j}(\partial_u^j R_i\cup \partial_s^j R_i).$$

Let $(Y,f)$ be a finitely presented system and $(\Sigma,\sigma)$ an SFT extension of $Y$ with  a factor map $\pi$ such that the image of each cylinder set for $i\in\mathcal{A}(\Sigma)$ is a rectangle $R_i$ in $Y$.  Following Fried in~\cite{Fri} we define the {\bf symbol set}, {\bf star}, and {\bf second star} for a set $Y'\subset Y$ in
the following manner:
$$
\begin{array}{llll}
\mathcal{A}(Y')=\{i\in\mathcal{A}(\Sigma)\, | \textrm{ there
exists a point }x\in Y'
\textrm{ where }x\in R_i\},\\
\mathrm{star}(Y')=\bigcup_{s\in\mathcal{A}(Y')}R_s,\textrm{ and}\\
\mathrm{star}_2(Y')=\bigcup_{y\in\mathrm{star}(Y')}\mathrm{star}(y).
\end{array}
$$

\section{Preliminary Results}

In this section we provide some useful results concerning factor maps from Smale spaces to finitely presented systems.  These results will be needed throughout the paper.  The next lemma provides a criteria on periodic points for a factor map from a Smale space to a finitely presented system to be infinite-to-one.  

\begin{lem}\label{keylemma}
Let $(X,F)$ be a transitive Smale space, $(Y,f)$ be finitely presented, and $\pi: X\rightarrow Y$ a factor map.  
If there exist convergent sequences $q_n\in Y$ and periodic points $p_n,p'_n\in X$ such that $\pi(p_n)=\pi(p'_n)=q_n$,   
$p_n\neq p'_n$ for all $n\in\mathbb{N}$, and $\lim_{n\rightarrow\infty}p_n=\lim_{n\rightarrow\infty}p'_n$, then 
$\pi$ is infinite-to-one.
\end{lem}

\noindent{\bf Proof.}  The idea of the proof is to use the POTP to construct infinitely many points in $X$ mapping to a point in $Y$.  Let $c_X$ and $c_Y$ be expansive constants for $X$ and $Y$, respectively.  Fix $\gamma>0$ such that 
$\gamma<c_X$ and if $x,x'\in X$ with $d(x,x')<\gamma$, then $d(\pi(x),\pi(x'))<c_Y/4$.  We now use the pseudo-orbit tracing property of $X$ and fix $\delta, \epsilon>0$ such that $\max\{\epsilon,\delta\}<\gamma/4$ and any $\epsilon$-pseudo orbit in $X$ is 
$\delta$-shadowed by an orbit in $X$.  We know there exists an $N\in\mathbb{N}$ such that $d(p_N,p'_N)<\epsilon$.  Then there 
exists an $M\in\mathbb{N}$ such that $$F^M(p_N)=p_N\textrm{ and }F^M(p'_N)=p'_N.$$  Furthermore, since $X$ is expansive there exists $$0< j<M 
\textrm{ such that }d(F^j(p_N),F^j(p'_N))>c_X.$$  Let $s, s'$ denote the finite sequences $$p_N, F(p_N),...,F^{M-1}(p_N)\textrm{ and }
p'_N, F(p'_N),...,F^{M-1}(p'_N),$$ respectively.  Let $\omega,\omega'$ be bi-infinite pseudo orbits that are concatenations of the segments $s$ and $s'$ where $\omega\neq\omega'$.  Then $\omega$ and 
$\omega'$ are 
$\delta$-shadowed by points $x$ and $x'$ in $X$, respectively.  Then 
$$\begin{array}{llll}
d(F^j(x), F^j(x')) & \geq d(F^j(p_N, f^j(p'_N))-2\delta\\
 &> c_X-\gamma/2>0
 \end{array}$$
and therefore $x\neq x'$.  Furthermore, $\pi(x)=\pi(x')=q_N$ by the choice of $\gamma$.  We know there are uncountably many such pseudo-orbits $\omega$, and therefore $\#(\pi^{-1}(q_N))=\infty$. $\Box$

We now use the previous lemma to prove a result about preimages of periodic points for finite-to-one maps between a Smale space and a finitely presented system.

\begin{lem}\label{lem5.2}
If $(X, F)$ is a transitive Smale space, $(Y,f)$ is finitely presented, $\pi$ is a finite-to-one factor map from $X$ to $Y$, and 
$M=\min_{y\in Y}\#(\pi^{-1}(y))$, then there exists a dense open set $W\subset Y$ such that each periodic point in $W$ has 
$M$ preimages.
\end{lem}

\noindent{\bf Proof.}  Let $y\in Y$ such that $\#\pi^{-1}(y)=M$.  Suppose there does not exist a neighborhood $U$ of $y$ such that 
each periodic point in $U$ has exactly $M$ preimages.  
Then there exists a sequence $q_n\rightarrow y$ such that each $q_n\in\mathrm{Per}(Y)$ has 
more than $M$ preimages.  
Since $\pi$ is continuous and $X$ compact we know that if $\{p_n\}\subset X$ is a sequence where $p_n\in \pi^{-1}(q_n)$ then any convergence subsequence of $\{p_n\}$ converges to a point of $\pi^{-1}(y)$.

Then there exists a point $x\in\pi^{-1}(y)$ and sequences $q_{n_j}\in Y$, $p_{n_j},p'_{n_j}\in X$ such that
\begin{itemize}
\item $q_{n_j}$ is a subsequence of $q_n$, 
\item $\pi(p_{n_j})=\pi(p'_{n_j})=q_{n_j}$,  
\item $p_{n_j}\neq p'_{n_j}$ for all $j\in\mathbb{N}$, and 
\item $\lim_{j\rightarrow\infty}p_{n_j}=\lim_{j\rightarrow\infty}p'_{n_j}=x$.  
\end{itemize}
Lemma~\ref{keylemma} then implies that $\pi$ is infinite-to-one, a contradiction.  Hence, there exists a neighborhood $U$ of $y$ such that each periodic point in $U$ has exactly $M$ preimages.  Finally, the transitivity of $X$ implies that $Y$ is transitive and $W=\bigcup_{n\in\mathbb{Z}}f^n(U)$ is a dense open set in $Y$ such that each periodic point in $W$ has exactly $M$ preimages.  $\Box$ 

The next lemma shows that the bracket is preserved under $\pi$ as long as the points in $X$ are sufficiently close.

\begin{lem}\label{lem3.4}
Let $\pi: X\rightarrow Y$ be a factor map from a Smale space $(X,F)$ to a
finitely presented system $(Y,f)$. Then there exists a constant
$\delta_{\pi}>0$ such that if $x,x'\in X$ with $d(x,x')<\delta_{\pi}$,
then $\pi([x,x'])=[\pi(x),\pi(x')]$.
\end{lem}

\noindent{\bf Proof.}   We assume that we are using adapted metrics
on both $Y$ and $X$.  Let $c$ be an expansive constant for $Y$ and
$\epsilon\in(0,c/2)$.  So $W^s_{\epsilon}(y)\cap W^u_{\epsilon}(y')$ consists of at most one point for all $y,y'\in Y$.  Let $\delta_{\pi},\epsilon'>0$ be sufficiently small such
that we have the following:  
\begin{itemize}
\item If $x,x'\in X$ and
$d(x,x')<\epsilon'$, then $d(\pi(x),\pi(x'))<\epsilon$.  
\item There exists a $\delta_{\pi}$ such that if points $x,x'\in X$ and
$d(x,x')<\delta_{\pi}$, then $W_{\epsilon'}^u(x)\cap
W^s_{\epsilon'}(x')$ consists of one point in $X$. 
\end{itemize}
Such a
$\delta_{\pi}$ and $\epsilon'$ exist since  $X$ is
a Smale space.

It then follows that if $d(x,x')<\delta_{\pi}$ that $[x,x']$ exists.  We
know that $$d(F^{-n}([x,x']),F^{-n}(x))<\epsilon'\textrm{ and }
d(F^{n}([x,x']),F^{n}(x'))<\epsilon'$$ for all $n\in\mathbb{N}$.  Since the
map $\pi$ is a factor map we have
$$
\begin{array}{llll}
d(\pi F^{-n}([x,x']),\pi F^{-n}(x))=d(f^{-n}(\pi[x,x']),f^{-n}\pi(x))
<\epsilon\textrm{ and}\\
d(\pi F^{n}([x,x']),\pi F^{n}(x'))=d(f^{n}(\pi[x,x']),f^{n}\pi(x'))
<\epsilon\textrm{ for all }n\in\mathbb{N}.
\end{array}$$
This implies that
$$\pi([x,x'])\in W^u_{\epsilon}(\pi(x))\cap W^s_{\epsilon}(\pi(x'))=
[\pi(x),\pi(x')].$$  $\Box$

\begin{cla}\label{c.equal}
 Let $(Y,f)$ be expansive with expansive constant $c>0$ and $\epsilon<c/2$.  If $(\Sigma, \sigma)$ is an SFT extension of $(Y,f)$ where the image of each cylinder set is a rectangle, then 
 $$\pi(W^i_{\mathrm{loc}}(t))=W^i(\pi(t), R_{t_0})$$
 for each $t\in \Sigma$ and $i=u$ or $s$.
 \end{cla}
 
 \noindent{\bf Proof.}  Let $\hat{t}\in W^u_{\mathrm{loc}}(t)$.  So $\hat{t}_j=t_j$ for all $j\leq 0$.  Hence, $$\pi(\sigma^{j}\hat{t})\in R_{t_j}\subset B_{\epsilon}(\pi(\sigma^jt))$$
for all $j\leq 0$.  This implies that $\pi(\hat{t})\in W^u_{\epsilon}(\pi(t))\cap R_{t_0}=W^u(\pi(t), R_{t_0})$.

Now suppose that $y\in W^u(\pi(t), R_{t_0})$.  So there exists $\hat{t}$ such that $\pi(\hat{t})=y$ and $\hat{t}_0=t_0$.  Let $s=[t,\hat{t}]$.  Then using a similar argument as in the previous lemma it is not hard to show that $\pi(s)=\pi([\hat{t}, t])=[\pi(\hat{t}), \pi(t)]=y$ and $s\in W^u_{\mathrm{loc}}(t)$.
$\Box$

\section{Finitely presented systems have resolving extensions to Smale spaces}

In this section we prove Theorem~\ref{thm1}.  The standard proof that every sofic shift  has a $u$-resolving extension to an SFT uses the future cover of the sofic shift, see~\cite[p. 75]{LM}.  
The idea of the proof of Theorem~\ref{thm1} is very similar to the future cover construction.  Namely, we will construct an SFT extension, related to the future cover, that codes the ``different'' unstable sets accumulating on a point.  The SFT cover is zero dimensional, hence, we will need to factor the SFT cover appropriately to a Smale space such that there is a $u$-resolving factor map from the Smale space to the finitely presented system.

More specifically, the proof of Theorem~\ref{thm1} will proceed in the following steps:  Let $(Y,f)$ be a finitely presented system with a one-step SFT, $(\Sigma, \sigma)$, an extension of $(Y,f)$,
such that the image of each cylinder set is a sufficiently small
rectangle. 
In each rectangle we know there is a product structure,  i.e., for a rectangle $R$ if $x,y\in R$, then $[x,y]$ and $[y,x]$ is defined.  Hence, the breakdown of a local product structure for a finitely presented system occurs at the boundaries of the rectangles.  We use the rectangles then to detect where the product structure breaks down.  We do this by coding the ``different" unstable sets accumulating on a point, as given by the rectangles, and 
construct a one-step SFT, $(\Sigma_+,\sigma)$,  and a
factor map $\pi_+^0$ from $\Sigma_+$ to $Y$.  
The coding given by the SFT, $(\Sigma_+, \sigma)$, is too crude, in general, and so we  define an
equivalence relation $E_{\alpha}\subset\Sigma_+\times\Sigma_+$ and a factor map
$\alpha:\Sigma_+\rightarrow X=\Sigma_+/E_{\alpha}$.  The equivalence relation $E_{\alpha}$ is defined so that if, on a uniformly small scale, the unstable set detected by $\omega\in\Sigma_+$ and $\omega'\in\Sigma_+$ is the same, then we say $(\omega, \omega')\in E_{\alpha}$.
Finally, we
define a factor map $\pi_+$ such that $\pi_+^0=\pi_+\alpha$ and show
that $\pi_+$ is $u$-resolving and $X$ is a Smale
space under the canonically induced map.  So we will have the
following diagram.

$$
\xymatrix{
& \Sigma_+ \ar[dr]^{\alpha}\ar[dd]^{\pi^0_+} &\\
\Sigma\ar[dr]_{\pi} & & X\ar[dl]^{\pi_+}\\
 & Y & }
$$

\subsection{Construction of $\Sigma_+$.}
We now proceed with the construction.  Let $d(\cdot,\cdot)$ be an adapted metric for $Y$ under the action
of $f$.  We first need the rectangles in $Y$ to be sufficiently small.  Fix $c>0$ an expansive constant
for $Y$, and $0<\epsilon<c/2$.  Let $(\Sigma, \sigma)$ be a one-step SFT extension of $Y$ such that  each cylinder set is a rectangle,
$$\max_{y\in Y}(
\mathrm{diam}(\mathrm{star}_2(y)))<\epsilon/5<c/10,$$
and for each $y\in Y\cap R_j$, where $R_j$ is the image of the cylinder $C_j$ under $\pi$, we have
$f(W^u(y, R_j))\subset W^u_{\epsilon}(y)$.  We note that the existence of such a $\Sigma$ is guaranteed by passing to a higher block presentation, if needed.


We now proceed with the construction of $(\Sigma_+, \sigma)$ using $(\Sigma, \sigma)$.  For $s\in\Sigma$ define
$$
E_i(s)=\{j\in\mathcal{A}(\Sigma)\, |\, W^u(\pi(\sigma^{-i}s),R_{s_i})\cap R_j\neq\emptyset\}.
$$
We will use the sets $E_i(s)$ to detect the different unstable sets accumulating on a point in $Y$.

\begin{rem}\label{rem1} From the definition it follows that $E_i(s)=E_{i-j}(\sigma^j(s))$ for all $s\in\Sigma$ and all $i,j\in\mathbb{Z}$.
\end{rem}

We now define the one-step SFT, $(\Sigma_+,\sigma)$, as an extension of
$(Y,f)$. Let $(\Sigma_+,\sigma)$ be an SFT with symbol set
$$\mathcal{A}(\Sigma_+)=\{(v,j)\, |\,(E_i(s),s_i)=(v,j)\textrm{ for some } s\in\Sigma\textrm{ and
}i\in\mathbb{Z}\}$$ and define allowable transitions from
$$(v,j)\in\mathcal{A}(\Sigma_+)\textrm{ to }
(v',j')\in\mathcal{A}(\Sigma_+)$$ if there exists an $s\in\Sigma$
and $i\in\mathbb{Z}$ such that $$(v,j)=(E_i(s),s_i)\textrm{ and }
(v',j')=(E_{i+1}(s),s_{i+1}).$$

A point $\omega\in \Sigma_+$ is {\bf canonically associated} with $s\in\Sigma$ if
$$\omega_i=(E_i(s),s_i)\textrm{ for all }i\in\mathbb{Z}.$$  
For each $\omega\in \Sigma_+$ we know that if $\omega_j=(v_j, i_j)$, then $i_j$ to $i_{j+1}$ is an allowed transition in $\Sigma$ for all $j\in\mathbb{Z}$.  Furthermore, since $(\Sigma, \sigma)$ is a one-step SFT we know that if  $\omega\in \Sigma_+$ and $\omega_j=(v_j,i_j)$ for all $j\in\mathbb{Z}$, then there exists a unique $s\in S$ such that $s_j=i_j$ for all $j\in\mathbb{Z}$.
Hence, there is
a canonical map $\pi^0_+:\Sigma_+\rightarrow Y$ such that
$\pi^0_+(\omega)=\pi(s)$.  It is clear from the definition of
$\Sigma_+$ that this is a well-defined map. 


\subsection{Properties of $\Sigma_+$.}  We now state and prove properties of $\Sigma_+$ that will be needed in the proof of Theorem~\ref{thm1}.  The next lemma shows
that the set of canonically associated $\omega$ is dense in
$\Sigma_+$.  

\begin{lem}\label{chainlem}  Let $\omega\in\Sigma_+$,
$k,l\in\mathbb{Z}$ with $k\leq l$.  Then there exists a point
$s\in\Sigma$ such that $(E_i(s),s_i)=\omega_i$ for all $k\leq
i\leq l$.
\end{lem}

\noindent{\bf Proof.}  The proof proceeds by induction.  From Remark~\ref{rem1} we may assume that $k=0$.  Fix $\omega\in\Sigma_+$.  For $l=0$ and $l=1$ we know, by the definition of $\Sigma_+$, that there exists an $s\in\Sigma$ such that $$(E_0(s),s_0)=\omega_0\textrm{ and }(E_1(s),s_1)=\omega_1.$$   Now suppose there exists an $s\in\Sigma$ such that $$(E_i(s),s_i)=\omega_i\textrm{ for all }0\leq i\leq m$$ where $m\geq 1$.  Fix such an $s$ and fix $s'\in\Sigma$ such that $$(E_m(s'),s'_m)=\omega_m\textrm{ and }(E_{m+1}(s'),s'_{m+1})=\omega_{m+1}.$$

We now construct the desired point by taking the bracket of $\sigma^ms'$ and $\sigma^ms$. 
Let $$s''=\sigma^{-m}[\sigma^{m}s',\sigma^{m}s].$$  From Claim~\ref{c.equal} we know that $$\sigma^is''\in W^u(\pi(\sigma^is), R_{s_i})\textrm{ and }s''_i=s_i\textrm{ for all }i\leq m.$$
Hence, 
$$(E_{i}(s''),s''_{i})=(E_i(s),s_i)=\omega_i\textrm{ for all }0\leq i\leq m.$$
We now show that $$(E_{m+1}(s''),s''_{m+1})=\omega_{m+1}.$$  We know that
$s''_{m+1}=s'_{m+1}$ so we need to prove that $E_{m+1}(s'')=E_{m+1}(s')$.

\begin{figure}
\begin{center}

\psfrag{1}{$R_{s_m}$}
\psfrag{2}{$R_n$}
\psfrag{3}{$f^{-1}(R_{s'_{m+1}})$}
\psfrag{4}{$f^{-1}(R_j)$}
\psfrag{a}{$\pi(\sigma^ms)$}
\psfrag{b}{$\pi(\sigma^ms'')$}
\psfrag{c}{$z_1$}
\psfrag{e}{$f^{-1}(y)$}
\psfrag{d}{$\pi(\sigma^ms')$}
\psfrag{f}{$W^s(\pi(\sigma^ms'), R_{s_m})$}
\psfrag{g}{$W^u(\pi(\sigma^ms'), R_{s_m})$}
\psfrag{h}{$W^u(\pi(\sigma^ms), R_{s_m})$}

\includegraphics[height=2.3in]{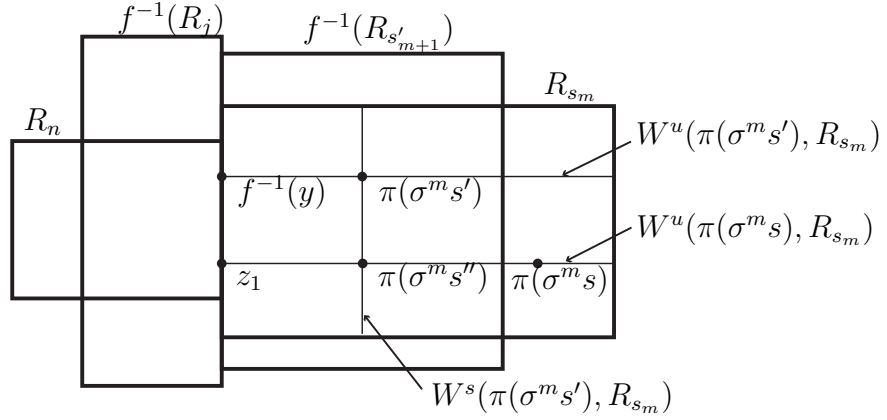}
 \caption{\label{f.equivalent}$j\in E_{m+1}(s'')$ }

\end{center}
\end{figure}

Fix $j\in E_{m+1}(s')$.  We now show that $j\in E_{m+1}(s'')$.  Let  
$$y\in W^u(\pi(\sigma^{m+1}s'), R_{s'_{m+1}})\cap R_j.$$  
We know from the Markov property that 
$$f^{-1}(W^u(y, R_j))\subset W^u(f^{-1}(y), R_n)$$ for some $n\in\mathcal{A}(\Sigma)$ where $n\rightarrow j$ is an allowed transition in $\Sigma$
and
$$f^{-1}(W^u(\pi(\sigma^{m+1}s'), R_{s'_{m+1}}))\subset W^u(\pi(\sigma^ms'), R_{s'_m}).$$
Hence, 
$$f^{-1}(y)\in W^u(\pi(\sigma^{m}s'), R_{s'_m})\cap R_n$$  
and $n\in E_m(s')=E_m(s)=E_m(s'')$.  

Hence, there exists a point 
$$
z_0\in W^u(\pi(\sigma^{m}s''), R_{s_m})\cap R_n
$$
and $$z_1=[f^{-1}(y), z_0]\in W^u(\pi(\sigma^{m}s), R_{s_m})\cap R_n,$$ see Figure~\ref{f.equivalent}.

Since 
$$z_1\in W^s(f^{-1}(y), R_n)\subset f^{-1}(W^s(y, R_j))$$ and
$$z_1\in W^s(f^{-1}(y), R_{s_m})\subset f^{-1}(W^s(y, R_{s'_{m+1}}))$$
we know that 
$$f(z_1)\in R_j\cap R_{s'_{m+1}}.$$
By construction we know that $f(z_1)\in W^u_{\epsilon}(\pi(\sigma^{m+1}s''))$.  Therefore, 
$$f(z_1)\in W^u(\pi(\sigma^{m+1}s''), R_{s'_{m+1}})\cap R_j$$
and $j\in E_{m+1}(s'')$.

A similar argument shows that $E_{m+1}(s')\subset E_{m+1}(s'')$.  Hence, $E_{m+1}(s'')=E_{m+1}(s')$. $\Box$

\begin{cla}\label{charactercla}
If $\omega\in \Sigma_+$, $\omega_0=(v,i)$, and $j\in v$, then there exists a point $y\in W^u(\pi^0_+(\omega),R_i)\cap R_j$.
\end{cla}

\noindent{\bf Proof.} The claim follows since canonically associated points are dense and the rectangles are closed.  Let $s_k\in\Sigma$ such that $\pi(s_k)\rightarrow \pi_+^0(\omega)$ and $\omega_k\rightarrow \omega$ where $\omega_k$ is the sequence in $\Sigma_+$ canonically associated with $s_k$ for all $k\in\mathbb{N}$.  We may assume $(\omega_k)_o=(v,i)$ for all $k\in\mathbb{N}$.  Then for each $k\in\mathbb{N}$ we know there exists $$y_k\in W^u(\pi(s_k),R_i)\cap R_j.$$  By possibly taking a subsequence we may assume that the $y_k$ are convergent to a point $y\in R_i\cap R_j$.  Since $R_i$ is a rectangle we know that $y\in W^u(\pi^0_+(\omega),R_i)\cap R_j$. $\Box$

From the definition of $\Sigma_+$ and the previous claim we know that the if  $(v,i)\in \mathcal{A}(\Sigma_+)$, then 
$$\pi^0_+(C_{(v,i)})=R_{(v,i)}$$ is a rectangle.  Furthermore,  
$R_{(v,i)}\subset R_i$ and if $y\in R_i\cap R_{(v,i)}$, then 
$$W^u(y,R_i)=W^u_{\epsilon}(y)\cap R_{(v,i)}=W^u(y,R_{(v,i)}).$$

To determine where the local product structure breaks down in $Y$ we want to extend the rectangles $R_{(v,i)}$.  To do this we define new rectangles, that we will denote $D_{(v,i)}$, that contain the sets $R_{(v,i)}$.  Let
$(v,i)\in\mathcal{A}(\Sigma_+)$ and define $$D_{(v,i)}=\bigcup_{z\in
R_{(v,i)}}(W^u_{\epsilon}(z)\cap(\bigcup_{j\in v} R_j)).$$

\begin{figure}[htb]
\begin{center}

\psfrag{1}{$R_{k}$}
\psfrag{2}{$R_i$}
\psfrag{3}{$R_j$}
\psfrag{a}{$y_2$}
\psfrag{b}{$z_2$}
\psfrag{c}{$z_3$}
\psfrag{d}{$y_3$}
\psfrag{e}{$y_1$}
\psfrag{f}{$z_1$}

\includegraphics[width=3in]{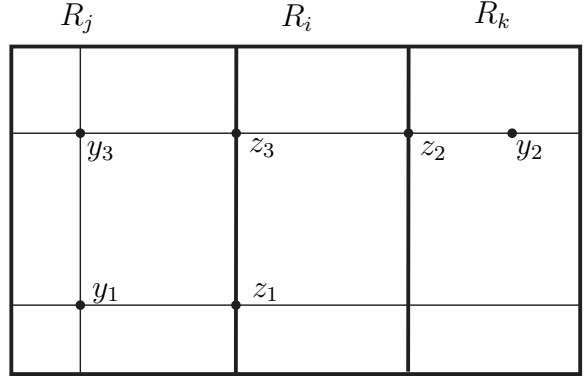}
\caption{$D_{(v,i)}$ is a rectangle} \label{figure1}
\end{center}
\end{figure}

\begin{cla}\label{rectanglecla}  For each $(v,i)\in\mathcal{A}(\Sigma_+)$
the set $D_{(v,i)}$ is a rectangle.  Furthermore, the rectangles have the Markov property.  
\end{cla}

\noindent{\bf Proof.}  The fact that the rectangles $D_{(v,i)}$ have the Markov property follows since the rectangles $R_i$ have the Markov property.
The set $D_{(v,i)}$ is a rectangle since canonically associated points are dense and the rectangles are closed.  More specifically, we know that for 
$y\in D_{(v,i)}$ we have
$$\mathrm{diam}(D_{(v,i)})\leq \mathrm{diam}(\mathrm{star}_2(y)) <
\epsilon.$$ Fix $y_1,y_2\in D_{(v,i)}$.  Then there exists $z_1,z_2\in
R_{(v,i)}$ and $j,k\in v$ such that $y_1\in W^u_{\epsilon}(z_1)\cap
R_j$ and $y_2\in W^u_{\epsilon}(z_2)\cap R_{k}$.  From
Claim~\ref{charactercla} it follows that there exists a point
$z_3\in W^u(z_2, R_i)\cap R_j$.  Then $[y_1,z_3]\in W^u_{\epsilon}(y_2)$, so 
$[y_1,z_3]=[y_1,y_2]=y_3$ exists and $$y_3\in W^u_{\epsilon}(z_2)\cap R_j$$ so $y_3\in D_{(v,i)}$, see Figure~\ref{figure1}.  Similarly, we can show that $[y_2,y_1]$ exists and is contained in $D_{(v,i)}$.  Hence, $D_{(v,i)}$ is a rectangle.
 $\Box$

\subsection{Construction of $X$.}  We now show that the rectangles $D_{(v,i)}$ allow us to define a uniform constant $\delta>0$ that can detect when two unstable sets coded by points in $\Sigma_+$ are different and when the unstable sets are simply coded differently.  We will use this constant to define the equivalence relation, $E_{\alpha}$, on points of $\Sigma_+$ such that $X=\Sigma_+/E_{\alpha}$.

Fix $\delta>0$ sufficiently small so that $2\delta<c$ and such that
the following hold:
\begin{enumerate}
\item  If $j,k\in\mathcal{A}(\Sigma)$ and $R_j\cap R_k=\emptyset$,
then $d(R_j,R_k)>2\delta$.
\item  If $y\in Y$ and $i\in\mathcal{A}(\Sigma)$ where $y\in
R_i$, then $W^u_{2\delta}(y)\subset\mathrm{star}(W^u(y,R_i))$.
\end{enumerate}
Property 2 says that if $s\in\Sigma$ and $y=\pi(s)$, then
$$W^u_{2\delta}(y)\subset W^u(y, D_{(E_0(s),s_0)}).$$
Define
$$B^u_{\delta}(\omega)=W^u_{\delta}(\pi^0_+(\omega))\cap
W^u(\pi^0_+(\omega), D_{\omega_0}).$$

\begin{cla}\label{usefulcla}
If $\omega\in\Sigma_+$, then $B^u_{\delta}(\omega)\subset
f(B^u_{\delta}(\sigma^{-1}\omega))$.
\end{cla}

\noindent{\bf Proof.}  We know from Claim~\ref{rectanglecla} that
$$f(W^u(\pi^0_+(\sigma^{-1}\omega),D_{\omega_{-1}})\supset
W^u(\pi^0_+(\omega),D_{\omega_0})$$ for all $\omega\in\Sigma_+$.
Since $Y$ has an adapted metric it also follows that
$$f(W^u_{\delta}(\pi^0_+(\sigma^{-1}\omega))\supset
W^u_{\delta}(\pi^0_+(\omega)).$$  Hence,
$B^u_{\delta}(\omega)\subset
f(B^u_{\delta}(\sigma^{-1}(\omega)))$.  $\Box$

\subsubsection{The relation $E_{\alpha}$}  We now define the relation $E_{\alpha}$ on $\Sigma_+$. 
\begin{defn} For
$\omega,\omega'\in\Sigma_+$ we say $(\omega,\omega')\in E_{\alpha}$, or
$\omega \sim_{\alpha}\omega'$, if $B^u_{\delta}(\sigma^j\omega)=B^u_{\delta}(\sigma^j\omega')$ for all  $j\in\mathbb{Z}$. 
\end{defn}
 Since $(Y, f)$ is expansive and
$\delta<c/2$ we have $$\pi_+^0(\omega)=\pi^0_+(\omega') \textrm{ for all }
\omega \sim_{\alpha}\omega'.$$  Hence, the relation $E_{\alpha}$ is
a subset of $E_{\pi^0_+}$. By definition the relation $E_{\alpha}$ is an equivalence relation.  The idea is that if $(\omega, \omega')\in E_{\alpha}$, then $\omega$ and $\omega'$ represent the same unstable sets at $\pi^0_+(\omega)$ and are simply coded differently.

\subsubsection{The space $X$}
Let $X=\Sigma_+/E_{\alpha}$ and $\alpha$ be the canonical map from
$\Sigma_+$ to $X$.  Furthermore, let the canonically induced action on $X$ be denoted by $F$.  Define $\pi_+$ to be the factor map from $X$ to $Y$ such
that $\pi^0_+=\pi_+\alpha$. We will use the metric on $X$ induced
from the metric on $\Sigma_+$ and the equivalence relation
$E_{\alpha}$. 

\begin{rem} If $Y$ is a sofic shift, then  the space $X$ is the Fischer cover as described in~\cite{Fisch2} and~\cite{Fisch1}.  
\end{rem}

\subsubsection{The map $\pi_+$ is $u$-resolving.}  The following definition and lemma
will be useful in proving $\pi_+$ is $u$-resolving.

\begin{defn}  An SFT $E'$ a subset of an SFT $E$ is {\bf forward closed}
if whenever $x\in E$ is backward asymptotic to $E'$, then $x\in
E'$.
\end{defn}

\begin{lem}~\cite{Boy1}\label{forwardclosed}
If $\gamma=\beta\alpha$ is a factorization of factor maps, then
$\beta$ is $u$-resolving if and only if $E_{\alpha}$ is forward
closed in $E_{\gamma}$.
\end{lem}

The next proposition will show that $\pi_+$ is $u$-resolving and $X$ is finitely presented.

\begin{prop}\label{closedprop} The relation $E_{\alpha}$ is closed, forward closed,
and a 1-step SFT.
\end{prop}

\noindent{\bf Proof.} To see that $E_{\alpha}$ is closed the idea is to take a sequence of points in $E_{\alpha}$ converging to some point in $\Sigma_+\times \Sigma_+$.
Since the unstable sets are all the same for the sequence this will also hold for the convergent point.  Indeed,
suppose
there is a sequence $(\omega_k,\omega'_k)\in E_{\alpha}$
converging to a point $(\omega,\omega')\in\Sigma_+\times\Sigma_+$.
Since $E_{\pi^0_+}$ is closed it follows that $(\omega,\omega')\in
E_{\pi^0_+}$ and $\pi^0_+(\omega)=\pi^0_+(\omega')$.

Fix a point $z\in B^u_{\delta}(\sigma^n\omega)$ for some
$n\in\mathbb{Z}$. Note that if $\omega_k$ converges to $\omega$
that $B^u_{\delta}(\sigma^n\omega_k)$ converges to
$B^u_{\delta}(\sigma^n\omega)$. It then follows that there
exists a sequence $z_k\in B^u_{\delta} (\sigma^n\omega_k)$
converging to $z$.  Since $B^u_{\delta}(\sigma^n\omega_k) =
B^u_{\delta}(\sigma^n\omega'_k)$ it follows that $z\in
B^u_{\delta}(\sigma^n\omega')$. Hence,
$$B^u_{\delta}(\sigma^n\omega)=B^u_{\delta}(\sigma^n\omega')
\textrm{ for all }n\in\mathbb{Z},$$
$(\omega,\omega')\in
E_{\alpha}$, and $E_{\alpha}$ is closed.

We now show that the equivalence relation $E_{\alpha}$ is forward
closed in $E_{\pi^0_+}$.  The idea is the following:  Suppose that $(\omega,\omega')\notin E_{\alpha}$ is backward asymptotic to $E_{\alpha}$.  Then the local unstable sets given by $B^u_{\delta}(\cdot)$ do not agree for some iterate of $\omega$ and $\omega'$.  Under backward iteration the point at which they do not agree will converge to the backward iterates of $\pi^0_+(\omega)$ and $\pi^0_+(\omega')$.  Hence,  any point in $E_{\alpha}$ for which $(\omega,\omega')$ is backward asymptotic will have local unstable sets given by $B^u_{\delta}(\cdot)$ will also differ, a contradiction.

More precisely,
let $(\omega,\omega')\in E_{\pi^0_+}$ be
backward asymptotic to $E_{\alpha}$.  Suppose that
$(\omega,\omega')\notin E_{\alpha}$.  Then there is a
$k\in\mathbb{Z}$ such that $B^u_{\delta}(\sigma^k\omega)\neq
B^u_{\delta}(\sigma^k\omega')$.    Since $B^u_{\delta}(\sigma^k\omega)\neq
B^u_{\delta}(\sigma^k\omega')$ we may assume there exists a $z\in Y$ such
that $$z\in B^u_{\delta}(\sigma^k\omega)\textrm{ and }z\notin
B^u_{\delta}(\sigma^k\omega').$$

From Claim~\ref{usefulcla} we know that
$$f^{-1}(B^u_{\delta}(\sigma^k\omega))\subset
B^u_{\delta}(\sigma^{k-1}\omega)\textrm{ and
}f^{-1}(B^u_{\delta}(\sigma^k\omega'))\subset
B^u_{\delta}(\sigma^{k-1}\omega').$$ This implies that
$$f^{-1}(z)\in
B^u_{\delta}(\sigma^{k-1}\omega).$$ Furthermore, we know that
$$W^u_{\delta}(\pi^0_+(\sigma^{k}\omega'))\subset
f(W^u_{\delta}(\pi^0_+(\sigma^{k-1}\omega')))\textrm{ and }$$
$$W^u(\pi^0_+(\sigma^k\omega'), D_{\omega'_k})\subset
f(W^u(\pi^0_+(\sigma^{k-1}\omega'), D_{\omega'_{k-1}})).$$ This
implies that
$$f^{-1}(z)\notin
B^u_{\delta}(\sigma^{k-1}\omega').$$  Continuing inductively it
follows that  $$f^{-n}(z)\in
B^u_{\delta}(\sigma^{k-n}\omega)\textrm{ and }f^{-n}(z)\notin
B^u_{\delta}(\sigma^{k-n}\omega')\textrm{ for all }n\in\mathbb{N}.$$
Since $z\in W^u_{\epsilon}(\pi^0_+(\sigma^{k}\omega))$ for all
$n\in\mathbb{N}$ it follows that $$d(f^{-n}(z),\pi^0_+(\sigma^{k-n}\omega))\rightarrow 0\textrm{ as }n\rightarrow \infty.$$ Fix
$J\in\mathbb{N}$ such that $f^{-n}(z)\in
B_{\delta/10}(\pi^0_+(\sigma^{k-n}\omega))$ for all $n\geq J$. Since $(\omega,\omega')$ is backward asymptotic to $E_{\alpha}$ we know there exists a point
$(\bar{\omega},\bar{\omega}')\in E_{\alpha}$ and a
subsequence of $\sigma^{-n}(\omega,\omega')$ converging to
$(\bar{\omega},\bar{\omega}')$.  Denote
$\bar{\omega}_0=(\bar{v},\bar{i})$ and
$\bar{\omega}'_0=(\bar{v}',\bar{i}')$.  Then for some $m\geq J$
sufficiently large we have
\begin{itemize}
\item $\omega_{k-{m}}=\bar{\omega}_0$,
\item $\omega'_{k-{m}}=\bar{\omega}'_0$,
\item $[\pi^0_+(\bar{\omega}),f^{-{m}}(z)]\in B^u_{\delta}(\bar{\omega})$,
and  
\item $[\pi^0_+(\bar{\omega}'),f^{-{m}}(z)]\notin
B^u_{\delta}(\bar{\omega}')$.  
\end{itemize}
This follows since
$D_{(\bar{v},\bar{i})}$ and $D_{(\bar{v}',\bar{i}')}$ are rectangles
and contradicts the fact that $(\bar{\omega},\bar{\omega}')\in
E_{\alpha}$.  Hence, $(\omega,\omega')\in E_{\alpha}$ and
$E_{\alpha}$ is forward closed in $E_{\pi^0_+}$.

We now show the set $E_{\alpha}$ is a 1-step SFT.  Let
$(a,b),(c,d),(e,f)\in\mathcal{A}(\Sigma_+)$ such that there exist
$(s,t),(s',t')\in E_{\alpha}$ where

$$\left.
\begin{array}{lllll}
 (a,b)& = & (s_{-1},t_{-1}) \\
(c,d) & = & (s_0,t_0)=  (s'_0,t'_0)\\
(e,f)&= &(s'_1,t'_1)  
\end{array}\right. .$$
We now show that $[(s,t),(s',t')]\in E_{\alpha}$.  Since
$E_{\pi^0_+}$ is a 1-step SFT it follows that $[(s,t),(s',t')]\in
E_{\pi^0_+}$.  Also $[(s,t),(s',t')]$ is backward asymptotic to
$E_{\alpha}$ so $[(s,t),(s',t')]\in E_{\alpha}$. Hence,
$(a,b)\rightarrow (c,d)\rightarrow (e,f)$ is an allowed transition
in $E_{\alpha}$ and $E_{\alpha}$ is a 1-step SFT. $\Box$

\subsubsection{Properties of $(X,F)$.}  
Before proceeding to the proof of Theorem~\ref{thm1} we show some of the properties of the space $(X,F)$. 

\begin{prop}
If $X=\Sigma_+/E_{\alpha}$, then $(X,F)$  is finitely presented and $\pi_+$ is $u$-resolving.  Furthermore, if
$s,s'\in\Sigma$ where $\pi(s)=\pi(s')=y$, then the points
$\omega,\omega'\in\Sigma_+$ canonically associated with $s$ and
$s'$, respectively, have $\alpha(\omega)=\alpha(\omega')$. 
\end{prop}

The first part of the proof of the above proposition follows from the previous result and Lemma~\ref{quotientlem}.  The second part of the above proposition follows from the definition of $E_{\alpha}$.

\begin{rem}Since the set of canonically associated points is dense in $\Sigma_+$ the above proposition implies that 
there is a canonical section $X'\subset X$ that is dense in $X$ and a bijective map $\tau:
Y\rightarrow X'$.  
\end{rem}

 We define rectangles in $X$ as the image of the cylinder sets in $\Sigma_+$.  
For $(v,i)\in\mathcal{A}(\Sigma_+)$ we let 
$$R_{(v,i)}=\{x\in X\,|\,\omega_0=(v,i),\alpha(\omega)=x\textrm{ for some }\omega\in\Sigma_+ \}$$
and if $x\in R_{(v,i)}$, then
$$W^u(x, R_{(v,i)})=\{ x'\in X\, |\, x'\in W^u(x), \pi_+(x')\in W^u(\pi_+(x), R_{(v,i)})\}.$$
The fact that these sets are rectangles follows directly from the definition.

The next lemma will be useful in showing that $(X,F)$
is a Smale space.  It will show that the rectangles $D_{(v,i)}$ are ``reconnected" by the map $\alpha$.

\begin{lem}\label{lem4.9} If $\pi_+(x)=\pi^0_+(\omega)=y$, $\alpha(\omega)=x$ and
$\omega_0=(v,i)$, then
$\pi_+(W^u(x,R_{\omega_0}))=W^u(y,R_i)$.  Furthermore,
for each $m\in v$ there exists a 
$(v_m,m)\in\mathcal{A}(\Sigma_+)$ such that $W^u(x,R_{\omega_0})\cap
R_{(v_m,m)}\neq\emptyset$.
\end{lem}

\noindent{\bf Proof.}  The first statement follows since 
$$
\begin{array}{llll}
W^u(y, R_i) &=W^u(\pi^0_+(\omega),R_{\omega_0})\\
& =\pi_+(W^u(\alpha(\omega), R_{\omega_0}))\\
& = \pi_+(W^u(x, R_{\omega_0})).
\end{array}
$$

We now show the second statement.  To do this we use points that are canonically associated.
For canonically associated points we will see that the definition of $E_{\alpha}$ implies the result.  Since the canonically associated points are dense we will see the result follows for all points in $X$.  Indeed,
let $s_k\in\Sigma$ be a sequence such that the sequence
$\omega_k\in\Sigma_+$ canonical associated to $s_k$ converges to
$\omega$.  Fix $m\in v$ and a sequence $s'_k\in
W^u_{\mathrm{loc}}(s_k)$ such that for each $k\in\mathbb{N}$ the
point $\pi(s'_k)\in R_i\cap R_m$.  By reducing to a subsequence, if
necessary, we may assume that the sequence $s'_k$ is convergent to a
point $s'\in\Sigma$ where $\pi(s')=y'\in W^u(y,R_i)\cap R_m$.

Let $\omega'_k\in\Sigma_+$ be the sequence canonically associated
with $s'_k$.  Then $\omega'_k\in W^u_{\mathrm{loc}}(\omega_k)$ for
all $k\in\mathbb{N}$ and $\omega'_k$ converges to some point
$\omega'\in W^u_{\mathrm{loc}}(\omega)$.  Let $x'=\alpha(\omega')$.
We want to show that there exists a
$(v_m,m)\in\mathcal{A}(\Sigma_+)$ such that $x'\in R_{(v_m,m)}$.

Fix a sequence $s''_k\in\Sigma$ such that $(s''_k)_0=m$ and
$\pi(s''_k)=\pi(s'_k)$ for all $k\in\mathbb{N}$.  Let
$\omega''_k\in\Sigma_+$ be the sequence canonically associated with
$s''_k$.  Again by possibly taking a subsequence we may assume that
$s''_k$ converges to a point $s''\in\Sigma$  and $\omega''_k$
converges to a point $\omega''$ where
$\pi(s'')=y'=\pi^0_{+}(\omega'')$.  Furthermore, there exists a
$(v_m,m)\in\mathcal{A}(\Sigma_+)$ such that $(\omega''_k)_0=(v_m,m)$
for all $k\in\mathbb{N}$. From the definition of
$\alpha$ we know that $(\omega''_k,\omega'_k)\in E_{\alpha}$ for all
$k\in\mathbb{N}$. Hence, $(\omega'',\omega')\in E_{\alpha}$ and $x'\in
R_{(v_m,m)}\cap W^u(x, R_{\omega_0})$. $\Box$

Let $x\in X$ and $\omega\in\Sigma_+$ where $x=\alpha(\omega)$ and
$\omega_0=(v,i)$.   For each $m\in v$ let $x_m\in W^u(x,
R_{\omega_0})\cap R_{(v_m,m)}$ and define
$$W^u(x,D_{\omega_0})=\bigcup_{m\in v}W^u(x_m,R_{(v_m,m)}).$$  From
the previous claim and lemma it is clear that $W^u(x,D_{\omega_0})$
is well defined and that the following proposition holds.

\begin{prop}\label{prop4.10} If $y=\pi_+(x)=\pi^0_{+}(\omega)$ and $\alpha(\omega)=x$, then $$W^u(y, D_{\omega_0})=\pi_+(W^u(x, D_{\omega_0})).$$
\end{prop}

Before proceeding to the proof of Theorem~\ref{thm1} we state a lemma that shows the space $X$ distinguishes the ``different'' unstable sets accumulating on a point of $Y$.

\begin{lem}~\label{l.character}
If $x_1,x_2\in X$ such that $\pi_+(x_1)=\pi_+(x_2)=y$, then $x_1\neq x_2$ if and only if $\pi_+(W^u(x_1))\neq\pi_+(W^u(x_2))$.
\end{lem}

\noindent{\bf Proof.}  Fix $\delta_0>0$ such that if $x,x'\in X$ and $d(x,x')<\delta_0$, then $$d(\pi_+(x),\pi_+(x'))<\delta.$$  Suppose there exist points $x_1$ and $x_2$ such that $\pi_+(W^u(x_1))\neq \pi_+(W^u(x_2))$.  Since $\pi_+$ is injective on $W^u(x_1)$ and $W^u(x_2)$ we may suppose, by possibly taking a backward iterate of $x_1$ and $x_2$, that there exists $$y'\in W^u_{\delta}(y)\cap \pi_+(W^u(x_1))\textrm{ and }y'\notin \pi_+(W^u(x_2)).$$  Let $x'_1=\pi_+^{-1}(y')\cap W^u(x_1)$.  By taking backwards iterates if needed we may assume that $x'_1\in W^u_{\delta_0}(x_1)$.  Let $\omega_1, \omega_2\in\Sigma_+$ such that $\alpha(\omega_1)=x_1$ and $\alpha(\omega_2)=x_2$.  Then 
$$
B^u_{\delta}(\omega_1)\supset \pi_+(W^u_{\delta_0}(x_1))\textrm{ and }
B^u _{\delta}(\omega_1)\neq B^u_{\delta}(\omega_2)$$ 
since $y'\in B^u _{\delta}(\omega_1)$ and 
$y'\notin B^u _{\delta}(\omega_2)$.  Hence, $x_1\neq x_2$ since $(\omega_1,\omega_2)\notin E_{\alpha}$.

Now suppose that $x_1\neq x_2$.  Then by taking backwards iterates, if necessary, we may assume that there exist $\omega_1,\omega_2\in\Sigma_+$ such that $\alpha(\omega_1)=x_1$, $\alpha(\omega_2)=x_2$, $B^u_{\delta}(\omega_1)\neq B^u_{\delta}(\omega_2)$, and there exists a point $y'\in W^u_{\delta}(y)$ such that 
$$y'\in B^u_{\delta}(\omega_1)\textrm{ and  }
y'\notin B^u_{\delta}(\omega_2).$$  
Hence, $\pi_+(W^u(x_1))\neq \pi_+(W^u(x_2))$ since $\pi_+^{-1}(B^u_{\delta}(\omega_1))$ and $\pi_+^{-1}(B^u_{\delta}(\omega_2))$ are neighborhoods of $x_1$ in $W^u(x_1)$ and $x_2$ in $W^u(x_2)$, respectively.  $\Box$

\subsection{Proof of Theorem~\ref{thm1}.}  Suppose $(X,F)$ is not
a Smale space.  Then there exists some $\epsilon'>0$ sufficiently
small and sequences $x_k$ and $x'_k$ in $X$ such that
$d(x_k,x'_k)<1/k$ and $W^s_{\epsilon'}(x_k)\cap
W^u_{\epsilon'}(x'_k)=\emptyset$.  By possibly taking a subsequence
we may assume that $x_k$ and $x'_k$ converge to a point $x\in X$,
and that $x_k\in R_{(v,i)}$ and $x'_k\in R_{(v',i')}$ for each
$k\in\mathbb{N}$ where $(v,i), (v',i')\in\mathcal{A}(\Sigma_+)$.

\begin{figure}[htb]
\begin{center}

\psfrag{1}{$R_{(v',i')}$}
\psfrag{2}{$R_{(v,i)}$}
\psfrag{a}{$x$}
\psfrag{b}{$x'_k$}
\psfrag{c}{$x_k$}
\psfrag{d}{$z'_k$}
\psfrag{e}{$z_k$}

\includegraphics[width=4in]{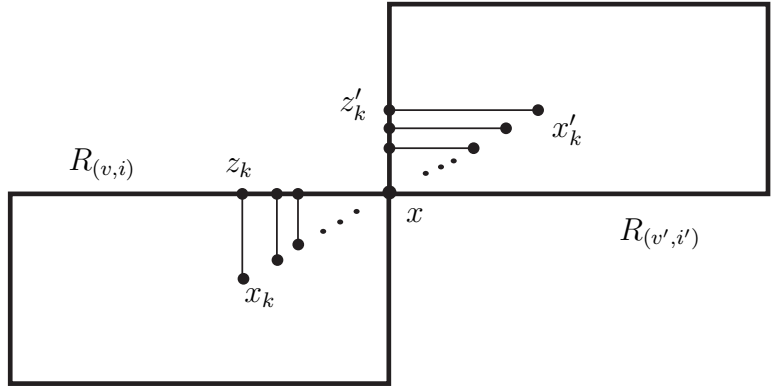}
\caption{$D_{(v,i)}$ is a rectangle} \label{l.product}
\end{center}
\end{figure}

By possibly taking a subsequence we may assume that the sequences
$z_k= W^u_{\epsilon'/2}(x)\cap W^s_{\epsilon'/2}(x_k)$ and $z'_k=
W^u_{\epsilon'/2}(x'_k)\cap W^s_{\epsilon'/2}(x)$ satisfy the
following, see Figure~\ref{l.product}:
\begin{itemize}
\item the sequences are
well defined,
\item $W^u_{\epsilon'/2}(z'_k)\cap
W^s_{\epsilon'/2}(z_k)=\emptyset$ for all $k\in\mathbb{N}$,
\item $z_k$ converges to $x$, and
\item $z'_k$ converges to $x$.
\end{itemize}

Let $\omega$ and $\omega'$ be in $\Sigma_+$ such that
$\alpha(\omega)=\alpha(\omega')=x$, $\omega_0=(v,i)$, and
$\omega'_0=(v',i')$.  From the above it follows that for all
$\epsilon''>0$ that $$W^u(x, D_{\omega_0})\cap B_{\epsilon''}(x)
\neq W^u(x, D_{\omega'_0})\cap B_{\epsilon''}(x).$$ The previous
proposition implies that $B^u_{\delta}(\omega)\neq
B^u_{\delta}(\omega')$, a contradiction.  Therefore, $(X,F)$ is a
Smale space. 

We now show that if $(Y,f)$ is transitive, then $(X,F)$ can be chosen to be transitive 
and $\pi_+$ a $u$-resolving one-to-one almost everywhere map.    
Since $(X,F)$ is a Smale space there is a
spectral decomposition such that $R(f)$ is decomposed into sets
$\Lambda_1,...,\Lambda_n$ such that $R(f)=\bigcup_{1\leq
i\leq n}\Lambda_i$ each set $\Lambda_i$ is closed, transitive, and
has dense periodic points.  Furthermore, we have
$$\mathrm{h}_{\mathrm{top}}(X)=\max_{1\leq i\leq
n}\mathrm{h}_{\mathrm{top}}(\Lambda_i).$$  Therefore, $X$ contains
an irreducible component $X_+$ of maximal entropy.  (For
definitions and properties of topological entropy see~\cite[p.
36]{BS}.)  From the Spectral Decomposition Theorem we know that $X_+$ is a transitive 
Smale space.  Furthermore, the image of $X_{+}$ under $\pi_{+}$ is $Y$ since $X_+$ is of maximal entropy.  

The next step is to show that $\pi_+$ is one-to-one on a residual set of $Y$.  For each $(v,j)\in\mathcal{A}(\Sigma_+)$ the set $R_{(v,j)}$ forms a closed rectangle and the set of these rectangles forms a 
cover $\mathcal{M}$ covering $Y$.  
The cover $\mathcal{M}$ is not necessarily a Markov partition since the rectangles may not be proper.  However, Fried 
in~\cite[p. 496-8]{Fri} shows there is an open and dense set $\mathcal{D}$ of $Y$ such that if $(v,j)\in\mathcal{A}(\Sigma_+)$ and  
$y\in \mathcal{D}\cap R_{(v,i)}$,  then $y\in\mathrm{int}(R_{(v,i)})$.  Let $W$ be an open and dense set in $Y$ such that each periodic point in $W$ has the minimal number of preimages for $\pi$.  If each point $p\in \mathrm{Per}(Y)\cap W$ has a unique pre-image under $\pi^0_+$, then it will follow that each point in $W\cap \mathcal{D}$ is contained in one rectangle $R_{(v,i)}$.  Hence, $V=\bigcap_{k\in\mathbb{Z}}f^k(\mathcal{D}\cap W)$ is a residual set in $Y$ and each point in $V$ would have a unique preimage in 
$X$ 
under $\pi_+$.  Hence, $X_+$ would be unique and each point in $V$ would have a unique preimage under $\pi_+$ in $X_+$.

To complete the proof we show that each point $p\in \mathrm{Per}(Y)\cap W$ has a unique pre-image under $\pi^0_+$.  Fix $p\in\mathrm{Per}(Y)\cap W$, let $s=\pi^{-1}(p)$, and $\omega\in\Sigma_+$ be the point canonically associated with $s$.  Suppose $\omega'\in \Sigma_+$ and $\pi^0_+(\omega')=p$.  We will show that $\omega'_0=\omega_0$.  Since $p$ was arbitrary we then have that $f^n(p)=\pi(\sigma^ns)$ for all $n\in\mathbb{Z}$ and $\sigma^n\omega$ is canonically associated with $\sigma^ns$.  Hence, $\omega'_n=\omega_n$ for all $n\in\mathbb{Z}$ and $\omega'=\omega$.

Let $\omega_0=(v,i)$ and suppose there exists some $\omega'\in\Sigma_+$ such that $p=\pi^0_+(\omega')$ and $\omega'_0=(v',i')$.  Then $i=i'$ since $p$ is contained in the interior of $R_i$.  Fix a sequence $s^k\in\Sigma$ such that the canonically associated sequence $\omega^k\in\Sigma_+$ converges to $\omega'$.  We may assume the $s^k$ are convergent, hence the sequence converges to $s$. 

We now show that $v=v'$.
Fix $j\in v'$.  Then there exists a sequence $t^k\in\Sigma$ such that 
$$\pi(t^k)\in W^u(\pi(s^k),R_i)\cap R_j.$$  
We may assume that $t^k$ converges to some point $t$.  Hence, 
$$\pi(t)\in W^u(\pi(s),R_i)\cap R_j$$ 
and $j\in v$.  

Now fix $j\in v$.  We will see that $j\in v'$ since $p$ is a periodic point contained in the interior of $R_i$ and the rectangles have the Markov property.  Indeed, let 
$$y\in W^u(p,R_i)\cap R_j$$
 and $N\in\mathbb{N}$ be the period of $p$.  Then $f^{-mN}(y)\in \mathrm{int}(R_i)$ for some $m\in\mathbb{N}$ since $p\in\mathrm{int}(R_i)$.  Let $t^k$ be a sequence of points in $\Sigma$ such that the canonically associated points $\omega^k$ converge to $\omega'$ and $t^k_{-lN}=i$ for $0\leq l\leq m$ and for all $k\in\mathbb{N}$.  Let $y_k=[y, \pi(t^k)]$.  From the Markov property of the rectangles we know that 
 $$f^{-mN}(y_k)=[f^{-mN}y, \pi(\sigma^{-mN}t^k)]$$ for all $k\in\mathbb{N}$.   Furthermore, by the Markov property of the rectangles and the fact that the rectangles are proper we know that 
 $$\begin{array}{llll}f^{-mN}(W^u(y, R_j))\subset W^u(p, R_i)\textrm{ and }\\ W^s(f^{-mN}(y), R_i)\subset f^{-mN}(W^s(y, R_j)).
 \end{array}$$
 Therefore, $f^{-mN}(y_k)\in f^{-mN}(W^s(y, R_j))$ and $y_k\in R_j$ for all $k\in\mathbb{N}$.  Hence, $j\in v'$ and $v=v'$ and $p$ has a unique preimage under $\pi^0_+$.
 
Let $\beta$ be the restriction of $\pi_{+}$ to $X_+$.  Since $\pi_{+}$ is $u$-resolving we know that 
$\beta$ is $u$-resolving.  Hence, $(X_+,\beta)$ is a transitive Smale space and $\beta$ is a $u$-resolving one-to-one 
almost everywhere factor map from $X_+$ to $Y$.  $\Box$

\section{Minimal resolving extensions for transitive finitely presented systems}

 In this section we prove Theorem~\ref{thm1a}.  Specifically, we show that if $(Y,f)$ is a transitive finitely presented system, then there exists a minimal $u$-resolving Smale space extension.\\

In general, a $u$-resolving factor map from a Smale space to a finitely presented system is injective, but not necessarily surjective, on unstable sets.  However, this is not the case for resolving maps between Smale spaces.  

\begin{defn} A factor map $\pi$ from $(X,F)$ to $(Y,f)$ is {\bf $u$-bijective} ({\bf $s$-bijective}) if for any $x\in X$ the map $\pi$ is a bijection from $W^u(x)$ ($W^s(x)$) to $W^u(\pi(x))$ ($W^s(\pi(x))$.
\end{defn}

\begin{thm}\label{t.bijection}\cite{Put3} If $(X,F)$ is a transitive Smale space and $(Y,f)$ is a Smale space, then a factor map $\pi$ is $u$-resolving ($s$-resolving) if and only if $\pi$ is $u$-bijective ($s$-bijective).
\end{thm}

We remark that it can be shown that if $\pi: X\rightarrow Y$ is a $u$-resolving factor map
from a transitive Smale space $(X,F)$ to a finitely presented system
$(Y,f)$, then there exists an open and dense set $W$in $Y$ such that if
$p\in \mathrm{per}(X)$ and $\pi(p)\in W$, then 
$\pi(W^u(p))=W^u(\pi(p))$.
\\

\noindent{\bf Proof of Theorem~\ref{thm1a}.}  Let $(Y_+,f_+)$ be the transitive Smale 
space extension of $(Y,f)$ constructed in the proof of Theorem~\ref{thm1} and 
$\pi_+$ the one-to-one almost everywhere factor map from $Y_+$ to $Y$.
Let $(X,g)$ be a transitive Smale space and $\alpha:X\rightarrow Y$ be a
$u$-resolving factor map.  

We will construct a map $\beta:
X\rightarrow Y_+$ such that $\beta$ is a $u$-resolving factor map.

$$
\xymatrix{
X\ar[dr]_{\alpha}\ar[rr]^{\beta} & & Y_+\ar[dl]^{\pi_+}\\
 & Y & }
$$
To construct the map $\beta$ we will use an irreducible component of the fiber product.

We now form the fibered product of $\pi_+:Y_+\rightarrow Y$ and $\alpha: X\rightarrow Y$.  Let
$$G=\{(x,z)\in X\times Y_+\, |\, \alpha(x)=\pi_+(z)\},$$
the map $p_1:G\rightarrow X$ be the coordinate projection onto $X$, and the map
$p_2:G\rightarrow Y_+$ be the coordinate projection onto $Y_+$.  Endow $G$ with the product metric.  Then the map $h:G\rightarrow 
G$ defined as $h(x,z)=(g(x),f_+(z))$ is a homeomorphism.  We then have the following commutative diagram.

$$
\xymatrix{
& G \ar[dl]_{p_1}\ar[dr]^{p_2} &\\
X\ar[dr]_{\alpha} & & Y_+\ar[dl]^{\pi_+}\\
 & Y & }
$$

It is not hard to see that $(G,h)$ is a Smale space and $p_1$ and $p_2$ are $u$-resolving.  Let $G_+$ be an irreducible component of maximal entropy.  We denote the projection from $G_+$ to $X$ as $\rho_1$ and the 
projection from $G_+$ to $Y_+$ as $\rho_2$.  Since the maps $\pi_+$ and $\alpha$ are 
$u$-resolving it follows that the maps $\rho_1$ and $\rho_2$ are $u$-resolving.  We will 
show that the map $\rho_1$ is injective, hence a homeomorphism.  The map $\beta=\rho_2\rho_1^{-1}$ 
will then be a $u$-resolving factor map from $X$ to $Y_+$.

We now proceed with the proof that $\rho_1$ is a homeomorphism. 
Since $X$, $Y_+$, and $G_+$ are transitive Smale spaces we know from Theorem~\ref{t.bijection} that 
$p_1(W^u(x,y))=W^u(x)$ and $p_2(W^u(x,y))=W^u(y)$ for any $(x,y)\in G_+$.
Suppose there exist points $(x_1,y_1),(x_2,y_2)\in G_+$ with the same image under $\rho_1$.  Then 
$$x_1=\rho_1(x_1,y_1)=\rho_1(x_2,y_2)=x_2.$$
Furthermore,  we have 
$$
\begin{array}{llll}
\pi_+(W^u(y_2)) & = \pi_+ \rho_2(W^u(x_1,y_2)) = \alpha \rho_1(W^u(x_1,y_2)) &\\
& =\alpha \rho_1(W^u(x_1,y_1)) = \pi_+ \rho_2(W^u(x_1,y_1) & = \pi_+(W^u(y_1)).
\end{array}$$
From Lemma~\ref{l.character} we know $y_1=y_2$.  Hence, 
$\rho_1$ is injective.  $\Box$

\section{Resolving extensions for almost Smale systems}

The following result gives another characterization of almost Smale systems.

\begin{prop}
Let $(X, F)$ be a transitive Smale space and $(Y,f)$ be a finitely presented
system.  A factor map $\theta: X\rightarrow Y$ is 1-1 on an open set if and
only if $\theta$ is $u$-resolving, $s$-resolving, and 1-1 on a residual
set.
\end{prop}

\noindent{\bf Proof.}  Suppose $\theta$ is 1-1 on an open set.  Then there
exist open sets $U\subset X$ and $V\subset Y$ such that $\theta(U)=V$ and
$\theta^{-1}(V)=U$.  Since, $f$ and $F$ are homeomorphisms it follows that
for all $i\in\mathbb{Z}$ that $\theta(F^i(U))=f^i(V)$ and
$\theta^{-1}(f^i(V))=F^i(U)$.  Since $(Y,f)$ is transitive it follows that
there is an open dense set $\bigcup_{i\in\mathbb{Z}}f^i(V)$ where $\theta$ is
1-1.

Choose $\epsilon>0$ and $\delta>0$ such that for all points $x,x'\in X$ where
$d(x,x')<\delta$ that $W^s_{\epsilon}(x)\cap W^u_{\epsilon}(x')$ consists of
one point.
Suppose $\theta$ is not $u$-resolving.  Then there exist $x_1, x_2\in X$ such that $x_1\in W^u(x_2)$ and
$\theta(x_1)=\theta(x_2)$.  We may assume that $d(x_1,x_2)<\delta/2$.  Let
$z$ be a transitive point for $(X,F)$ such that the $d(z,x_1)<\delta/2$.
It follows that the points $[x_1,z]=z_1$ and $[x_2,z]=z_2$ map to the
same point under $\theta$.

Let $x\in U$ and $r>0$ such that $B_{4r}(x)\subset U$.  Choose
$n\in\mathbb{N}$ so that 
\begin{itemize}
\item $d(F^{-n}(z),x)<r$, 
\item $d(F^{-n}(z),F^{-n}(z_1))<r$,
and  
\item $d(F^{-n}(z),F^{-n}(z_1))<r$.
\end{itemize}  
Then the points $F^{-n}(z_1),F^{-n}(z_2)
\in U$ and $\theta(F^{-n}(z_1))=\theta(F^{-n}(z_2))$, a contradiction.
Hence, $\theta$ is $u$-resolving.  Similarly, one can show $\theta$ is
$s$-resolving.

Suppose $\theta$ is
$u$-resolving, $s$-resolving, and 1-1 on a residual set.  Let $z$ be a
transitive point in $X$ and let $\epsilon>0$ and $\delta>0$ be given by
the canonical coordinates on $X$.  Define the rectangle
$$R_{\delta}(z)=\{ x\in X\,\,| \, x=[z_1,z_2]\textrm{ for some }
z_1\in W^u_{\delta/2}(z)\textrm{ and }z_2\in W^s_{\delta/2}(z)\}.$$
The rectangle $R_{\delta}(z)$ is a proper rectangle.  Suppose that
there exist points $x_1,x_2\in\mathrm{int}(R_{\delta}(z))$ such that
$\theta(x_1)=\theta(x_2)$.  Then the points
$[x_1,z], [x_2,z]\in W^u_{\delta/2}(z)$ get mapped to the same point
under $\theta$ and $\theta$ is not $u$-resolving, a contradiction.
Therefore, $\theta$ is 1-1 on $\mathrm{int}(R_{\delta}(z))$.  $\Box$

\begin{prop}\label{constanttoone}
Let $(X, f)$ and $(Y, g)$ be transitive Smale spaces and $\pi:X\rightarrow
Y$ be a $u$-resolving and $s$-resolving factor map.  Then $\pi$ is
constant to one.
\end{prop}

\noindent{\bf Proof.}  The map $\pi$ is bounded to one from a
result in~\cite{Boy1}.  Let $y\in Y$ such that $\#(\pi^{-1}(y))=k$
is maximal and fix $x\in\pi^{-1}(y)$.  Fix $\epsilon>0$ and $\delta>0$ from
the canonical coordinates on $X$.  We may suppose that $\delta<\epsilon$.

We now show that $\pi$ is 1-1 on $B_{\delta/2}(x)$.  Let $x'\in X$ and
$x_1, x_2\in B_{\delta/2}(x')$
such that $\pi(x_1)=\pi(x_2)$.  This implies that
$\pi([x_1,x'])=\pi([x_2,x'])$.  Since $\pi$ is $u$-resolving we have
$[x_1,x']=[x_2,x']$, which implies that $x_1\in W^s_{\epsilon}(x_2)$.  Since
$\pi$ is $s$-resolving we then have $x_1=x_2$ and $\pi$ is 1-1 on $B_{\delta/2}(x)$.

Let $V=\bigcup_{x\in\pi^{-1}(b)}B_{\delta/2}(x)$.  Then each point
in $V$ has exactly $k$ pre-images for $\pi$.  Since $Y$ is
transitive it follows that there is an open and dense set of
points in $Y$ with $k$ pre-images.  Hence any transitive point has
$k$ pre-images.

Fix $y\in Y$, a transitive point $z\in Y$, and sequence $n_i$ such
that $g^{n_i}(z)$ converges to $y$.  Each $g^{n_i}(z)$ has $k$
pre-images which we can order $z^1_i,..., z^k_i$.  By compactness
of $X$ we may suppose, by possibly taking a subsequence of $g^{n_i}(z)$, that each
of the sequences $z^j_i$ converge to a point in $\pi^{-1}(y)$. Since
each $\delta/2$ neighborhood of a point in $\pi^{-1}(y)$ is 1-1,
as shown in the previous paragraph, it follows that $y$ has $k$
pre-images.  Hence, $\pi$ is constant to one. $\Box$\\

\noindent{\bf Proof of Theorem~\ref{thm2}.}  Let $F$ be the fibered product
for $\theta$ and $\phi$.  Denote $\rho_1$ and $\rho_2$ the canonical
projections for $F$ onto the first and second coordinates, respectively.

$$
\xymatrix{
& F \ar[dl]_{\rho_1}\ar[dr]^{\rho_2} &\\
X'\ar[dr]_{\phi} & & X\ar[dl]^{\theta}\\
 & Y & }
$$

Since $\theta$ is $u$-resolving, $s$-resolving, and 1-1 on an open set
it follows that
$\rho_1$ is.  From Proposition~\ref{constanttoone} it follows that $\rho_1$
is $k$ to one for some constant $k\in\mathbb{N}$.
Since $\theta$ is one to
one on an open set it follows that $k=1$.  Therefore, $\rho_1$ is invertible.
If $\rho=\rho_2\rho_1^{-1}$, then $\phi=\theta\rho$.  The map $\rho$ is clearly
continuous, commuting, and onto.  Hence, $\rho$ is a factor map and $\phi$
factors through $\rho$.  $\Box$

\section{Lifting factor maps through 
u-resolving maps to s-resolving maps}

We now show that every finite-to-one factor map between transitive finitely presented systems lifts through $u$-resolving maps  to an $s$-resolving map between Smale spaces.
The proof of Theorem~\ref{thm2a} will proceed in two steps.  First, we prove the following proposition.

\begin{prop}\label{diagprop}
Let $(X,f)$ be a transitive finitely presented system, $(Y,g)$ be a finitely presented system and $\pi$ be a finite-to-one factor map from $X$ to $Y$.  Then there exist transitive Smale
spaces $(\bar{X}, \bar{f})$ and $(\bar{Y}, \bar{g})$ and finite-to-one factor maps $\gamma$, $\beta$, and $\bar{\pi}$
such that the following diagram commutes.
$$
\xymatrix{
\bar{X} \ar[d]_{\gamma}\ar[r]^{\bar{\pi}} & \bar{Y}\ar[d]^{\beta}\\
X\ar[r]_{\pi} &  Y} $$
Moreover, the maps $\gamma$ and $\beta$ are $u$-resolving.
\end{prop}

\noindent{\bf Proof.}  Since $X$ is transitive under $f$ and $\pi$ is a factor map we know that $Y$ is transitive under $g$.
Let $(X', F)$ and $(\bar{Y}, \bar{g})$ be the minimal transitive Smale space extension of $(X, f)$ and $(Y, g)$, respectively, constructed in the proof of Theorem~\ref{thm1}.  Furthermore, let $\pi_+$ and $\beta$ be the $u$-resolving almost
one-to-one factor map from $X'$ to $X$ and $\bar{Y}$ to $Y$, respectively.  The map $\pi'=\pi\pi_+$ is a finite-to-one factor map from $X'$ to $Y$.  

We now form the fibered product of $\bar{Y}$ and $X'$ as in the last section.  Let
$$G=\{(x,z)\in X'\times\bar{Y}\, |\, \pi(x)=\beta(z)\},$$
$p_1:G\rightarrow X'$ the coordinate projection onto $X$, and
$p_2:G\rightarrow \bar{Y}$ the coordinate projection onto $\bar{Y}$.  
$$
\xymatrix{
G\ar[d]_{p_1}\ar[r]^{p_2} &  \bar{Y}\ar[dd]^{\beta}\\
X'\ar[d]_{\pi_+}\ar[dr]^{\pi'} & \\
X\ar[r]_{\pi} & Y}
$$
We endow $G$ with the product metric.  Then the map $h:G\rightarrow
G$ defined as $h(x,z)=(F(x),\bar{g}(z))$ is a homeomorphism.

As in the previous section it is not hard to see that $(G,h)$ is a Smale space and $p_1$ is $u$-resolving.  Let $\bar{X}$ be an irreducible component of maximal entropy, $\psi$ be the restriction $p_1$ to $\bar{X}$, and $\bar{\pi}$ be the restriction of $p_2$ to $\bar{X}$.   We then have the following commutative diagram.
$$
\xymatrix{
\bar{X}\ar[d]_{\psi}\ar[r]^{\bar{\pi}} &  \bar{Y}\ar[dd]^{\beta}\\
X'\ar[d]_{\pi_+}\ar[dr]^{\pi'} & \\
X\ar[r]_{\pi} & Y}
$$
Then $\gamma=\pi_+\psi$ is surjective and is a $u$-resolving factor map from $\bar{X}$ to $X$.  The map $\bar{\pi}$ is surjective and a finite-to-one factor map. 
$\Box$

Theorem~\ref{thm2a} will follow from the above result and the next theorem.

\begin{thm}\label{thm3}
Let $(\bar{X}, \bar{f})$ and $(\bar{Y}, \bar{g})$ be transitive Smale spaces and $\bar{\pi}$ a finite-to-one factor map from $\bar{X}$ to $\bar{Y}$.  Then there exists transitive Smale spaces $(\tilde{X}, \tilde{f})$ and $(\tilde{Y},\tilde{g})$ and finite-to-one factor maps $\tilde{\gamma}$, $\tilde{\beta}$, and $\tilde{\pi}$
such that the following diagram commutes.

$$
\xymatrix{
\tilde{X} \ar[d]_{\tilde{\gamma}}\ar[r]^{\tilde{\pi}} & \tilde{Y}\ar[d]^{\tilde{\beta}}\\
\bar{X}\ar[r]_{\bar{\pi}} &  \bar{Y}} $$
Moreover, the maps $\tilde{\gamma}$ and $\tilde{\beta}$ are $u$-resolving, and $\tilde{\pi}$ is $s$-resolving.
\end{thm}

Assuming the above result we now prove Theorem~\ref{thm2a}. \\

\noindent{\bf Proof of Theorem~\ref{thm2a}.}  From Theorem~\ref{thm3} and Proposition~\ref{diagprop} we have the following commutative diagram.
$$
\xymatrix{
\tilde{X} \ar[d]_{\tilde{\gamma}}\ar[r]^{\tilde{\pi}} & \tilde{Y}\ar[d]^{\tilde{\beta}}\\
\bar{X} \ar[d]_{\gamma}\ar[r]^{\bar{\pi}} & \bar{Y}\ar[d]^{\beta}\\
X\ar[r]_{\pi} &  Y} $$
Where $(\tilde{X}, \tilde{f})$, $(\tilde{Y}, \tilde{g})$, $(\bar{X}, F)$, and $(\bar{Y}, \bar{g})$ are transitive Smale spaces, and $X$ and $Y$ are transitive and finitely presented.  The maps $\pi$, $\tilde{\pi}$, $\gamma\tilde{\gamma}$, and $\beta\tilde{\beta}$ are finite-to-one factor maps.  Moreover, the map $\tilde{\pi}$ is $s$-resolving, and the maps $\gamma\tilde{\gamma}$ and $\beta\tilde{\beta}$ are $u$-resolving.
$\Box$

The proof of Theorem~\ref{thm3} proceeds in a few steps.  The first part is to build a $u$-resolving extension of $\bar{Y}$.  We will not use the minimal extension, but instead construct a different $u$-resolving extension.  The construction is similar to the construction in the proof of Theorem~\ref{thm1}, however, we define a different relation $E_{\theta}$ that is a subset of  $E_{\alpha}$.  We do this to help define the $s$-resolving map.  The next step is to form an appropriate fiber product.  In the last step we extend the notion of a magic word from symbolic dynamical systems to maps between finitely presented systems.  This concept of a magic word is then used to construct the appropriate commuting diagrams.\\

\noindent{\bf Proof of Theorem~\ref{thm3}.}  Let $c>0$ be an expansive constant for both $X$ and $Y$. Let $c/2>\epsilon\geq\delta_0>0$ such that for any points $x,x'\in\bar{X}$ and $y,y'\in\bar{Y}$ where $d(x,x')<\delta_0$ and $d(y,y')<\delta_0$ that $W^s_{\epsilon}(x)\cap W^u_{\epsilon}(x')$ and $W^s_{\epsilon}(y)\cap W^u_{\epsilon}(y')$ consists of one point in $\bar{X}$ and $\bar{Y}$, respectively.

Let $\Sigma$ be a transitive SFT such that there is a factor map $\nu:\Sigma\rightarrow \bar{X}$.  Furthermore, we may assume that $\Sigma$ generates a Markov partition on $\bar{X}$ and that the rectangles formed by the image of the cylinder sets in $\Sigma$ are sufficiently small so that for any $x\in\bar{X}$ and $y\in\bar{Y}$ we have $\mathrm{diam}(\mathrm{star}_2(x))<\delta_0/10$ and
$\mathrm{diam}(\mathrm{star}_2(y))<\delta_0/10$.

Define an SFT $(\Sigma_+,\sigma)$ an extension of $\bar{Y}$ with factor map $\pi^0_+$, and $\delta>0$ as in the proof of Theorem~\ref{thm1}.  For any $\omega\in\Sigma_+$ let $B^u_{\delta}(\sigma^j_+(\omega))$ be defined as in the proof of Theorem~\ref{thm1}.  We now define a relation $E_{\theta}$ on the points of $\Sigma_+$ such that we have the following commutative diagram.
$$
\xymatrix{
 & & \Sigma_+ \ar[dr]^{\theta}\ar[dd]_{\pi^0_+} &\\
 & & & Y_+\ar[dl]_{\pi_+}\\
\Sigma\ar[r]^{\nu} &  \bar{X}\ar[r]^{\bar{\pi}} & \bar{Y} & }
$$
The map $\pi_+$ will be $u$-resolving and $Y_+$ will be a Smale space.  However, in general, $Y_+$ will not be the minimal Smale space extension.  More specifically, we construct $Y_+$ such that if $t,t'\in \Sigma$ where $t\neq t'$ satisfy $\bar{\pi}\nu(t)=\bar{\pi}\nu(t')$ and $t\in W^s(t')$, then $\theta(\omega)\neq \theta(\omega')$ where $\omega$ and $\omega'$ are in $\Sigma_+$ and are the canonically associated
points of $t$ and $t'$, respectively.

Define
$$V^u_{\delta}(\sigma^j_+\omega)=\{k\in\mathcal{A}(\Sigma)\,|\,\omega_j=(v_j,i_j), k\in v_j, \textrm{ and }R_k\cap B^u_{\delta}(\sigma_+^j\omega)\neq\emptyset\}.$$
\begin{figure}

\psfrag{1}{$R_1$}
\psfrag{2}{$R_2$}
\psfrag{3}{$R_3$}
\psfrag{y}{$y$}
\psfrag{A}{$W^u(y, R_1)$}

\includegraphics[height=2.5in]{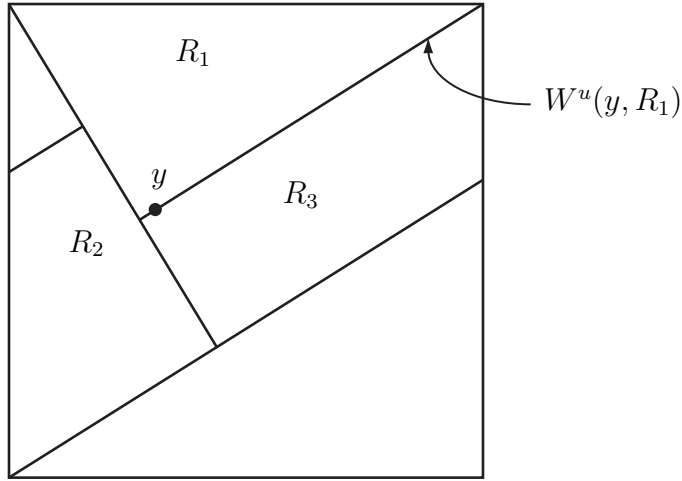}
 \caption{\label{f.theta}Preimages for $\theta$ }

\end{figure}
For any two points $\omega,\omega'\in\Sigma_+$ we say
$$
(\omega,\omega')\in E_{\theta}\textrm{ if }V^u_{\delta}(\sigma^j_+\omega)=V^u_{\delta}(\sigma^j_+\omega)
\textrm{ for all }j\in\mathbb{Z}.$$
By the definition of $E_{\theta}$ we see that $E_{\theta}$ is an equivalence relation.  For points $\omega, \omega'\in \Sigma_+$ where $(\omega,\omega')\in E_{\theta}$ we notice that it is not enough simply that $B^u_{\delta}(\sigma^j_+\omega)=B^u_{\delta}(\sigma^j_+\omega')$ for all $j\in\mathbb{Z}$, we need the sets to be in the same rectangles in $v_j$ and $v'_j$ where $\omega_j=(v_j, j)$ and $\omega'_j=(v'_j, j')$ for all $j\in\mathbb{Z}$.  

To show how this differs from the definition of $\alpha$ notice that if $f_a$ is hyperbolic toral automorphism of $\mathbb{T}^2$ and $\Sigma$ generates a sufficiently small Markov partition of $\mathbb{T}^2$ and $\Sigma_+$ and $\alpha$ are defined as in the proof of Theorem~\ref{thm1}, then we know that the set $\Sigma_+/E_{\alpha}$ and the induced action on this space is the minimal $u$-resolving extension to a Smale space.  Hence, $\Sigma_+/E_{\alpha}=\mathbb{T}^2$ , the induced action is just $f_a$, and $\pi_+$ is the identity map.  However, the space $\Sigma_+/E_{\theta}$ will not be $\mathbb{T}^2$.  
Specifically, Figure~\ref{f.theta} shows a neighborhood of a point $y$ in $\mathbb{T}^2$ and rectangles for a Markov partition for $f_a$.  Then the point $y$ will have at least three preimages in $\Sigma_+/E_{\theta}$.  Namely, one containing $1$ and $2$ in $V^u_{\delta}(\cdot)$ but not $3$, one containing $1$, $2$, and $3$ in $V^u_{\delta}(\cdot)$, and one containing $2$ and $3$ but not $1$ in $V^u_{\delta}(\cdot)$. 

Let $Y_+=\Sigma_+/E_{\theta}$ and $\theta$ be the canonical factor map from $\Sigma_+$ to $Y_+$.  Let $\pi_+$ be the canonical factor map so $\pi^0_+=\pi_+\theta$.  We will show that $E_{\theta}\subset E_{\alpha}$, the set $Y_+$ is a Smale space under the action of the induced map, and $\pi_+$ is $u$-resolving.  The proof of these facts is very similar to the proof in Theorem~\ref{thm1} that $\Sigma_+/E_{\alpha}$ is a Smale space and the resulting factor map is $u$-resolving.  

\begin{lem}
The map $\pi_+$ is $u$-resolving, and $Y_+$ is a Smale space.
\end{lem}

\noindent{\bf Proof.}   From the size of $\delta$ and $\mathrm{star}_2(y)$ for all $y\in \bar{Y}$ we know that if $(\omega,\tilde{\omega})\in E_{\theta}$, then $\pi_+^0(\omega)=\pi_+^0(\tilde{\omega})$.  Furthermore, if $V^u_{\delta}(\sigma^j_+\omega)=V^u_{\delta}(\sigma^j_+\tilde{\omega})$, then $B^u_{\delta}(\sigma^j_+\omega)=B^u_{\delta}(\sigma^j_+\tilde{\omega})$.  It then follows that $E_{\theta}\subset E_{\alpha}$.

We first show $E_{\theta}$ is closed.
Let $\omega^k,\tilde{\omega}^k$ be sequences in $\Sigma_+$ converging to $\omega$ and $\tilde{\omega}$, respectively, such that $(\omega^k,\tilde{\omega}^k)\in E_{\theta}$.  Then $(\omega,\tilde{\omega})\in E_{\alpha}$.  We may assume that $\omega^k_0=\omega_0$ and $\tilde{\omega}^k_0=\tilde{\omega}_0$ for all $k$.

Suppose $(\omega,\tilde{\omega})\notin E_{\theta}$.  We may assume that $V^u_{\delta}(\omega)\neq V^u_{\delta}(\tilde{\omega})$ and there exists a $$j\in V^u_{\delta}(\omega)\textrm{ such that }j\notin V^u_{\delta}(\tilde{\omega}).$$  Then for $k$ sufficiently large we have $B^u_{\delta}(\omega_k)\cap R_j\neq\emptyset$ since
\begin{itemize}
\item
$B^u_{\delta}(\omega_k)$ is open relative to $W^u(\pi^0_+(\omega_k), D_{(\omega_k)_0})$,
\item $R_j$ is closed, and
\item the unstable sets vary continuously inside $D_{(\omega_k)_0}$.
\end{itemize}
Hence, $j\in V^u_{\delta}(\omega_k)=V^u_{\delta}(\tilde{\omega}_k)$.  If $\tilde{\omega}_0=(\tilde{v},\tilde{i})$, then $j\in \tilde{v}$.  Since $\pi^0_+(\omega)=\pi^0_+(\tilde{\omega})$ and $B^u_{\delta}(\omega)\cap R_j\neq\emptyset$ we have $j\in V^u_{\delta}(\tilde{\omega})$.

We now show that $E_{\theta}$ is forward closed.  Let $(\omega,\omega')\in E_{\pi^0_+}$ be backward asymptotic to $E_{\theta}$.  Then $(\omega,\omega')\in E_{\alpha}$.  Suppose for some $j\in\mathbb{Z}$ that $V^u_{\delta}(\sigma^j\omega)\neq V^u_{\delta}(\sigma^j\omega')$.  If $V^u_{\delta}(\sigma^{j-1}\omega)= V^u_{\delta}(\sigma^{j-1}\omega')$, then using the Markov property of $\Sigma$ we know that $V^u_{\delta}(\sigma^j\omega)= V^u_{\delta}(\sigma^j\omega')$.  Hence, for all $n\in\mathbb{N}$ we have
$V^u_{\delta}(\sigma^{j-n}\omega)\neq V^u_{\delta}(\sigma^{j-n}\omega')$.  We may suppose there is a $$l_0\in V^u_{\delta}(\sigma^j\omega)\textrm{ such that }  l_0\notin V^u_{\delta}(\sigma^j\omega')$$ and $y_0\in R_{l_0}\cap B^u_{\delta}(\sigma^j\omega)$.  Then $y_1=f^{-1}(y_0)\in R_{l_1}$ such that $l_1$ to $l_0$ is an allowed transition for $\Sigma$ and $l_1\in V^u_{\delta}(\sigma^{j-1}\omega)$.  Furthermore, $l_1\notin V^u_{\delta}(\sigma^{j-1}\omega')$ since otherwise $l_0$ would be in $V^u_{\delta}(\sigma^j\omega')$, and
$$d(\pi^0_+(\sigma^{j-1}\omega), y_1)\leq \lambda d(\pi^0_+(\sigma^j\omega), y_0)$$
where $\lambda\in (0,1)$ is a fixed constant depending on the adapted metric $d(\cdot, \cdot)$.  Continuing inductively we see that for each $n\in\mathbb{N}$ there is a $l_n\in\mathcal{A}(\Sigma)$ such that $$l_n\in V^u_{\delta}(\sigma^{j-n}\omega)\textrm{ and }l_n\notin V^u_{\delta}(\sigma^{j-n}\omega'),$$
and a point $y_n=f^{-n}(y_0)$ such that $y_n\in R_{l_n}$ and
$$d(\pi^0_+(\sigma^{j-n}\omega), y_n)\leq \lambda^n d(\pi^0_+(\sigma^j\omega), y_0).$$

Fix a subsequence $(\sigma_+^{-n_j}\omega,\sigma_+^{-n_j}\omega')$ converging to a point $(\bar{\omega},\bar{\omega'})\in E_{\theta}$.  We may suppose that $l_{n_j}$ are all the same.  Since $$y_{n_j}\rightarrow \pi^0_+(\bar{\omega})= \pi^0_+(\bar{\omega'})\textrm{ as }j\rightarrow\infty$$ we know that $$l_{n_j}\in V^u_{\delta}(\bar{\omega})\textrm{ and }l_{n_j}\notin V^u_{\delta}(\bar{\omega'}).$$
Indeed his shows $\pi^0_+(\bar{\omega})\in R_{l_{n_j}}$ and if $\bar{\omega}_0=(\bar{v},\bar{i})$ and $\bar{\omega'}_0=(\bar{v'},\bar{i})$, then $l_{n_j}\in \bar{v}$ and $l_{n_j}\in \bar{v'}$, a contradiction.  Therefore, $(\omega,\omega')\in E_{\theta}$ and $E_{\theta}$ is forward closed.

From Lemma~\ref{forwardclosed} we know that $\pi_+$ is $u$-resolving and $Y_+$ is finitely presented.  We now show that $Y_+$ is a Smale space.  If $t,t'\in\Sigma$ such that $\bar{\pi}\nu(t)=\bar{\pi}\nu(t')$, and $\omega,\omega'\in\Sigma_+$ are canonically associated with $t$ and $t'$, respectively, then the definition of $V^u_{\delta}$ implies that $(\omega,\omega')\in E_{\theta}$.  Then we can prove similar statements to Lemma~\ref{lem4.9} and Proposition~\ref{prop4.10} replacing $\theta$ with $\alpha$.  Hence, using the same argument as in the proof of Theorem~\ref{thm1} we can show that $Y_+$ is a Smale space.
$\Box$

Now let $\tilde{Y}$ be a transitive component of maximal entropy in $Y_+$ from the Spectral Decomposition Theorem.  Let $\tilde{\beta}$ be the restriction of $\pi_+$ to $\tilde{Y}$.  As stated previously this implies that $\tilde{\beta}: \tilde{Y}\rightarrow \bar{Y}$ is a finite-to-one $u$-resolving factor map.  Let $\Sigma'_+$ be an irreducible component of $\Sigma_+$ such that $\theta(\Sigma'_+)=\tilde{Y}$.

Let $G_0$ be the fiber product of $\Sigma$ and $\tilde{Y}$ and maps $\rho_1$ and $\rho_2$ the projections onto $\Sigma$ and $\tilde{Y}$, respectively.  Fix $G$
a transitive component of $G_0$ of maximal entropy from the Spectral Decomposition Theorem and let $\pi_1$ and $\pi_2$ be the restriction of $\rho_1$ and $\rho_2$, respectively, to $G$.  
$$
\xymatrix{
 & & G\ar[dll]_{\pi_1}\ar[dr]^{\pi_2} &\\
\Sigma\ar[r]^{\pi}\ar[drr]_{\tilde{\pi}\pi} & \bar{X}\ar[dr]^{\bar{\pi}} & &\tilde{Y}\ar[dl]_{\tilde{\beta}}\\
 & & \bar{Y} &}
$$

We now construct a factor of $G$ that will give us the desired commuting diagram.  Let $\tau$ send a point $(s,y)\in G$ to $(\pi(s),y)$ and let $\tilde{X}$ be the image of $G$ under $\tau$.  It is not hard to see that the space $\tilde{X}$ is a transitive Smale space since $\bar{X}$ is a Smale space and a factor of $\Sigma$ by $\pi$.  Furthermore, there exist maps $\tilde{\gamma}$ and $\tilde{\pi}$ projections from $\tilde{X}$ onto $\bar{X}$ and $\tilde{Y}$, respectively.  We then have the following commutative diagram.

$$
\xymatrix{
 & & G\ar[ddll]_{\pi_1}\ar[ddr]^{\pi_2}\ar[d]_{\tau} &\\
 & & \tilde{X}\ar[dl]^{\tilde{\gamma}}\ar[dr]_{\tilde{\pi}} &\\
\Sigma\ar[r]^{\pi}\ar[drr]_{\tilde{\pi}\pi} & \bar{X}\ar[dr]^{\bar{\pi}} & &\tilde{Y}\ar[dl]_{\tilde{\beta}}\\
 & & \bar{Y} &}
$$

\begin{lem}\label{c.sresolving}
The map $\pi_2$ is $s$-resolving.
\end{lem}

We postpone the proof of the above lemma to the next subsection and
proceed with the proof of the theorem.
Since $\tilde{\beta}$ is $u$-resolving we know that $\pi_1$ and $\tilde{\gamma}$ are $u$-resolving.  The maps $\tau$ and $\tilde{\pi}$ are $s$-resolving since $\pi_2$ is $s$-resolving and $\tilde{\pi}\tau=\pi_2$.  From Theorem~\ref{t.bijection} we know for each $(t,y)\in G$ we have
$\tau (W^s(t,y))=W^s (\tau(t),y)$.  Therefore, $\tilde{\pi}$ is $s$-resolving and we have the desired commutative diagram. $\Box$

\subsection{Proof of Lemma~\ref{c.sresolving}}

Before proceeding to the proof of Lemma~\ref{c.sresolving} we extend the notion of a magic word in symbolic dynamics to finite-to-one maps from SFT's to Smale spaces.  For a definition and basic properties of a magic word for finite-to-one maps between shift spaces see~\cite[p. 303]{LM}.   These arguments will show that each transitive point $y$ in $Y$ under $f$ has a unique preimage for $\pi_+$and that $\bar{Y}$ is the unique irreducible component in $Y_+$ and the set $G$ is transitive so $G=G_0$.  From this we will be able to conclude that $\pi_2$ is $s$-resolving.

\begin{cla}
Let $\Sigma$ be a transitive one-step SFT and $\phi$ be a finite-to-one factor map from $\Sigma$ to a transitive expansive system $(Y,f)$.  Then there exists $K\in\mathbb{N}$ so that if $w$, $w'$, and $w''$ are words of size $K$, and if $w$ is pairwise related to $w'$ by $\pi$, and $w'$ is pairwise related to $w''$ by $\pi$, then $w$ is pairwise related to $w''$ by $\pi$.
\end{cla}

\noindent{\bf Proof.}  Suppose that no such $K$ exists.  Then for each $n\in\mathbb{N}$ sufficiently large there exists words $w_n$, $w'_n$, and $w''_n$ where $|w_n|\geq 2n$ and such that  $w_n$ is pairwise related to $w'_n$, and $w'_n$ is pairwise related to $w''_n$, and $w_n$ is not pairwise related to $w''_n$.  Then for each $n$ there exists points $s^n, (s')^n, (s'')^n\in \Sigma$ such that $s^n[-n,n]\subset w_n$, $(s')^n[-n,n]\subset w'_n$, and $(s'')^n[-n,n]\subset w''_n$.  Since $\Sigma$ is compact we know that the sequence $(s^n, (s')^n, (s'')^n)\in \Sigma^3$ has a convergent subsequence $(s^{n_j}, (s')^{n_j}, (s'')^{n_j})\in \Sigma^3$ converging to a point $(s, s', s'')\in\Sigma^3$.  Since $w_{n_j}$ and $w'_{n_j}$ are pairwise related for all $j\in\mathbb{N}$ and $w'_{n_j}$ and $w''_{n_j}$ are pairwise related for all $j\in\mathbb{N}$ we know that $\phi(s)=\phi(s')$ and $\phi(s')=\phi(s'')$.  So $\phi(s)=\phi(s'')$ and for all $j$ sufficiently large we have $w_{n_j}$ pairwise related to $w''_{n_j}$, a contradiction.  Hence, such a $K\in\mathbb{N}$ exists. $\Box$

Fix a $K\in\mathbb{N}$ as in the previous claim.   
Let $m,n\in\mathbb{Z}$ where $m\leq -K$ and $K\leq n$.  Let $w$ be a word such that there exists an $s\in\Sigma$ where $s[m,n]=w$ and let $\mathcal{W}_{m,n}(w)$ be the maximal family of words pairwise related words by $\phi$ containing $w$.  Let $m\leq j\leq n$ and
$$\mathrm{deg}(\mathcal{W}_{m,n}(w),j)=\#\{\bar{w}_j\,| \bar{w}\in\mathcal{W}_{m,n}(w)\}.$$  We define $$\mathrm{deg}(\mathcal{W}_{m,n}(w))=\min_{m\leq j\leq n}\mathrm{deg}(\mathcal{W}_{m,n}(w),j)$$ and  $d=\min_{\mathcal{W}_{m,n}(w)}\mathrm{deg}(\mathcal{W}_{m,n})$.
Let $D=\min\{\#\phi^{-1}(y)\,|\, y\in Y\}$.  Then we know that $D\geq d$.  We will show that $D=d$.

Fix $\mathcal{W}_{m,n}(w)$ such that  $d=\mathrm{deg}(\mathcal{W}_{m,n}(w))$.  By possibly decreasing $m$, increasing $n$, and shifting the $j$th coordinate of $w$ we may assume that $\mathrm{deg}(\mathcal{W}_{m,n}(w))=\mathrm{deg}(\mathcal{W}_{m,n}(w),0)$.
Then for words $\bar{w},\bar{w}'\in\mathcal{W}_{m,n}(w)$ and $u$ such that $\bar{w}u\bar{w}'$ is a word in $\Sigma$ we define $\mathcal{W}^*(u)$ to be the set of $\Sigma$ words $w_1u_1w'_1$ such that $u$ and $u_1$ are pairwise related and $w_1,w'_1\in\mathcal{W}_{m,n}(w)$.  Therefore,  by the size of $K$ we know that words in $\mathcal{W}^*(u)$ are pairwise related.  Fix a $u$ such that $\mathcal{W}^*(u)$ is nonempty.

\begin{cla} $\mathcal{W}^*(u)$ is a maximal family of related words.
\end{cla}

\noindent{\bf Proof.}  Suppose $x$ is a word pairwise related to a word $w^*=w_1u_1w'_1\in W^*(u)$.   Then $x=x_1y_1x'_1$ where $|x_1|=|x'_1|=|n-m|$ and $|y_1|=|u|$.  From the size of $K$ we know $x_1$ and $x'_1$ are pairwise related to $w_1$ and $w'_1$, respectively.  Hence, $x_1,x_2\in W_{m,n}(w)$.  This implies that $x\in W^*(u)$.   $\Box$

So $\mathrm{deg}(\mathcal{W}^*(u), 0)\geq d$ and $\mathrm{deg}(\mathcal{W}^*(u), 0)\leq \mathrm{deg}(\mathcal{W}_{m,n}, 0)=d$.  Therefore, $\mathrm{deg}(W^*(u), 0)=d$.

\begin{cla}\label{c.magic2} 
Let $\{i_1,...,i_d\}=\{ \bar{w}_0 \, | \bar{w}\in\mathcal{W}_{m.n}(w)\}$ and $w^*\in \mathcal{W}^*(u)$.  Then there is a permutation $\tau$ of $\{i_i,...,i_d\}$ such that for $1\leq j\leq d$ there exists a unique word $v=v(i_j)$ such that $i_j v \tau(i_j)=w^*[0,|u|+n+1]$.
\end{cla}

\noindent{\bf Proof.}  We first show that $v$ is unique.  Suppose for some $j\in\{1,...d\}$ there exists words $i_jv\tau(i_j)$ and $i_jv'\tau(i_j)$.  Then there exist
points $t,t'\in\Sigma$ such that
\begin{itemize}
\item $t\neq t'$,
\item $t[0,|u| +n +1]=i_jv\tau(i_j)$, 
\item $t'[0,|u| +n +1]=i_jv'\tau(i_j)$, and
\item $t_i=t'_i$ for $i<0$ and $i>|v| +2$.
\end{itemize}
Then $\phi(t)=\phi(t')$ and $t\in W^u(t')\cap W^s(t')$.  From Lemma 2.1 in~\cite{Put} we know $t=t'$, a contradiction.

We now show that $v$ exists.  For $1\leq j\leq d$ there exists some $w^*\in W^*(u)$ such that $w^*_0=i_j$.  Otherwise, as stated above $\mathrm{deg}(\mathcal{W}^*(u),0)<d$, a contradiction.  Similarly, there exists some $w^*\in W^*(u)$ such that $w^*_{|u|+1}=i_j$.

Suppose for some $1\leq j\leq d$ there exists $k_1,k_2\in\{1,...d\}$ and words $v'$ and $v''$ such that $i_j v' i_{k_1}=w_1^*[0,|u|+n +1]$ and $i_j v'' i_{k_2}=w_2^*[0,|u|+n +1]$ for some $w_1^*, w_2^*\in\mathcal{W}^*(u)$.

Define $\mathcal{X}$ to be the collection of all words $w^*[1,|u|+n+1]= vi_s$ for some $w^*\in\mathcal{W}^*(u)$.  Since $n\geq K$ we know that each of the words in $\mathcal{X}$ are pairwise related.  Let $\mathcal{X}^k$ be the collection of allowed concatenations of words $x_1\cdots x_k$ where each $x_j\in \mathcal{X}$.  Then any two elements in $\mathcal{X}^k$ are pairwise related for each $k\in\mathbb{N}$.  We claim that $\#(\mathcal{X}^{k})<\#(\mathcal{X}^{k+1})$.   Let $x^{(k)}\in\mathcal{X}^k$ and $x^{(k)}$ end in $i_s$.  Then there exists some word $w^*\in\mathcal{W}^*(u)$ where $w^*[0,|u|+n+1]=i_svi_t$ for some $t\in\{1,...,d\}$.  Therefore, we know that each word $x^{(k)}\in \mathcal{X}^k$ has an extension to a word $x^{(k+1)}\in \mathcal{X}^{k+1}$ and $\#(\mathcal{X}^{k})\leq \#(\mathcal{X}^{k+1})$.  Now 
let $x^{(k)}\in\mathcal{X}^k$ such that $x^{(k)}$ ends in $i_j$.   Then $$x^{(k)}v' i_{k_1},x^{(k)}v'' i_{k_2} \in\mathcal{X}^{k+1}$$ and are distinct.  Hence, $\#(\mathcal{X}^{k})<\#(\mathcal{X}^{k+1})$.  

This implies that $\mathcal{X}^{D+1}$ has more than $D$ elements all of which can be extended to points in $\Sigma$ mapping to the same point in $Y$,  a contradiction.
$\Box$

\begin{cla}
Choose $t\in\Sigma$ such that a word in $\mathcal{W}_{m,n}(w)$ occurs infinitely often in the past and future.  Then $\#(\phi^{-1}(\phi( t)))=d$ and there does not exist a point $t'\in\Sigma$ such that $t'\in W^s(t)\cap \phi^{-1}(\phi(t))$ and $t'\neq t$.
\end{cla}

\noindent{\bf Proof}  Let $\bar{w}$ be the word in $\mathcal{W}_{m,n}(w)$ occurring infinitely often in $t$.  We know that there exists $N,M\in\mathbb{N}$ such that $t[-N, M]=\bar{w}u\bar{w}$ for some word $u$.  From the previous claim there exist $d$ words related to $t[-N,M]$.  As $N$ and $M$ can be chosen arbitrarily large we know that $\phi^{-1}(\phi(t))=d$.

Now suppose that $t'\in W^s(t)\cap \phi^{-1}(\phi(t))$.  We may suppose that $t'_i=t_i$ for all $i\geq 0$ and $t[0,m+n]=\bar{w}$.  Then there exists $N\in\mathbb{N}$ such that $t[-N, m+n]=\bar{w}u\bar{w}$ for some word $u$ and $t'[-N, m+n]=w'u'\bar{w}$ where $w'\in\mathcal{W}_{m,n}(w)$.  From the previous claim we know that $t'[-N+m, m+n]=t[-N+n,m+n]$.  Since $N$ can be made arbitrarily large we have $t'=t$. $\Box$

In particular, the above claim shows that $$d\geq\min_{y\in Y}\#\phi^{-1}(y)=D$$ so $d=D$.
Furthermore, if $t' \in \phi^{-1}(\phi( t))$, then $t$ and $t'$ are pairwise related,  $t'$  has words of $\mathcal{W}_{m,n}$ occurring infinitely often in the past and future, and the words occur in the same location as those of $t$.

The next claim says that if there is a finite-to-one factor map from a transitive SFT to a finitely presented system and two points on the same stable set map to the same point, then for some iterate the points are on the unstable boundaries of their respective rectangles.  This will be a helpful characterization in the proof of Theorem~\ref{thm2a} since $\theta$ was defined to separate such points.  This will imply that $\pi_2$ is $s$-resolving.

\begin{cla}\label{c.magic}  Suppose $\Sigma$ is a transitive one-step SFT, $Y$ is a transitive Smale spaces, and $\phi:\Sigma\rightarrow Y$ is a finite-to-one factor map.  If $t,t'\in\Sigma$ such that $\phi(t)=\phi(t')=y$ and $t\in W^s(t')$, then for some $j\in\mathbb{Z}$ the point $$\phi(\sigma^j(t'))\in\partial^{t_j}_u R_{t'_j}\cap \partial ^{t'_j}_u R_{t_j}.$$
\end{cla}

\noindent{\bf Proof.}  We want to show that 
$$\phi(t)=\phi(t')\in\partial^{t_0}_u R_{t'_0}$$ so there exists a sequence 
$$y_i\in (W^s_{\epsilon}(\phi(t'))\cap R_{t_0})-R_{t'_0}$$ converging to $\phi(t')$.  The argument that 
$$\phi(t)=\phi(t')\in\partial^{t'_0}_u R_{t_0}$$
will be similar and left to the reader.

We may suppose that $t_i=t'_i$ for all $i\geq 1$ and $t_0\neq t'_0$.  Let $\mathcal{W}_{m,n}(w)$ have degree $d$.  Since $\Sigma$ is transitive we know that for each $M\in\mathbb{N}$ there exists a $t^M\in\Sigma$ such that $t^M_i=t_i$ for all $-M\leq i\leq M$ and words from  $\mathcal{W}_{m,n}(w)$ occur infinitely often in the past and in the future for $t^M$.


Let $s^M=[t,t^M]$.  Then $\phi(s^M)\in W^s(\phi(t), R_{t_0})$.  We will show that $$\phi(s^M)\notin W^s(\phi(t'),R_{t'_0}).$$  As $M\rightarrow\infty$ we know that $\phi(s^M)\rightarrow \phi(t)$.  This implies that $\phi(t)\in \partial^{t_0}_u R_{t'_0}$.

Suppose that $$\phi(s^M)\in W^s(\phi(t'),R_{t'_0}).$$  Then there exists $s\in\Sigma$ such that
\begin{itemize}
\item $s_0=t'_0$,
\item $\phi(s)=\phi(s^M)$, and
\item $\phi(s)\in W^s(\phi(t'), R_{t'_0})$.
\end{itemize}
Let $x[0,\infty)$ be the infinite sequence such that
$$x_0=t'_0,\textrm{ and } x_i=t^M_i\textrm{ for all }i\geq 1.$$
Let $s'\in\Sigma$ such that $s'[0,\infty)=x[0,\infty)$ and $s'_i=s_i$ for all $i\leq 0$.  We know that for all $i\leq 0$ we have $$s'_i\sim_{\phi}s^M_i=t^M_i,$$ and $s'_i=t^M_i$ for all $i\geq 1$.  Furthermore,
$$(s'_i,t^M_i)\rightarrow (s'_{i+1},t^M_{i+1})$$
is an allowed transition in $E_{\phi}$ for all $i\in\mathbb{Z}$.
Then $$\phi(s')=\phi(t^M)\textrm{ and }s'\in W^s(t^M),$$ a contradiction.  Hence, $\phi(s^M)\notin W^s(\phi(t'),R_{t'_0})$. $\Box$

Before proceeding to the proof of Lemma~\ref{c.sresolving} we need
to prove some additional properties concerning pre-images of
periodic points and transitive points under $\pi^0_+$.

\begin{cla}\label{c.unique}
Let $\Sigma$ be a transitive $1$-step SFT and
$\phi:\Sigma\rightarrow Y$ is a finite-to-one factor map onto a
finitely presented system $(Y,f)$. If $p\in\mathrm{Per}(Y,f)$ and
there exist $s,s'\in\Sigma$ such that $\phi(s)=\phi(s')=p$ and
$s_0=s'_0$, then $s=s'$.
\end{cla}

The proof of Claim~\ref{c.unique} is similar to the proof of Lemma~\ref{keylemma} and is left to the reader.  

\begin{cla}\label{c.a7}
Let $\Sigma$ be a transitive $1$-step SFT, $(Y,f)$ be finitely
presented, and $\phi:\Sigma\rightarrow Y$ be finite-to-one.  Let
$D$ be the minimum number of pre-images under $\phi$, the set $W\subset Y$ be open and dense such that each periodic point in $W$ has $D$ preimages under $\phi$ and every point in $W$ is contained in the interior of a rectangle.  If 
$p\in W$ is a periodic point, then $\#
(\pi^0_+)^{-1}(p)=D$.
\end{cla}

\noindent{\bf Proof.}  From Lemma~\ref{lem5.2} we know that the set $W$ exists.  We know that for any rectangle containing $p$ that $p$ is
contained in the interior of the rectangle.  Let
$\omega\in\Sigma_+$ such that $\pi^0_+(\omega)=p$ and let
$\omega_0=(v,i)$.  Then there exists some $s\in\phi^{-1}(p)$ such
that $s_0=i$.  Let $\omega'\in\Sigma_+$ be canonically associated
with $s$ and $\omega'_0=(v',i)$. Fix a sequence $s_k\in\Sigma$ such
that the canonically associated sequence $\omega_k\in\Sigma_+$
converges to $\omega$. We may assume the $s_k$ are convergent.  From the previous claim we know that the sequence $s_k$ converges to $s$.

Suppose there exists a $j\in v$ such that $j\notin v'$.  Then
there exists a sequence $t_k\in\Sigma$ such that $\phi(t_k)\in
W^u(\phi(s_k),R_i)\cap R_j$.  We may assume that $t_k$ converges
to some point $t$.  Hence, $\phi(t)\in W^u(\phi(s),R_i)\cap R_j$
and $j\in v'$, a contradiction.

Suppose there exists a $j\in v'$ such that $j\notin v$.  Let $y\in
W^u(p,R_i)\cap R_j$ and $N\in\mathbb{N}$ be the period of $p$.
Then $f^{-mN}(y)\in \mathrm{int}(R_i)$ for some $m\in\mathbb{N}$
since $p\in\mathrm{int}(R_i)$.  Fix such an $m\in\mathbb{N}$. Then
there exists a $K\in\mathbb{N}$ such that $(s_k)_{-lN}=i$ for all
$0\leq l\leq m$ where $k\geq K$.  Such a $K$ exists since
$\sigma^{-lN}(s_k)\rightarrow s$ as $k\rightarrow \infty$ for all
$l$.  For $k\geq K$ let
$$y_k=[f^{-mN}(y),\pi(\sigma^{-mN}(s_k))]\in R_i.$$  
Since $f^{-mN}(y)\in R_i\cap f^{-mN}(R_j)$ we know from the Markov property of the rectangles that $f^{mN}(y_k)\in R_j$ for all $k\geq K$.  Furthermore, 
we know that
$$f^{lN}(y_k)=[f^{(l-m)N}(y),\pi(\sigma^{(l-m)N}(s_k)]$$ for all
$0\leq l\leq m$ since this point is unique. Then $$f^{mN}(y_k)\in W^u(\pi(s_k),
R_i)\cap R_j$$ for all $k\geq K$. This implies that $j\in v'$, a
contradiction.  Therefore, $v=v'$ and $\omega_0=\omega'_0$.  Since $p$ is periodic we see by iterating $p$ that $\omega=\omega'$. $\Box$

\begin{rem}\label{r.magic}
The previous claim shows that if $\Sigma$ is a transitive $1$-step SFT, $(Y,f)$ is finitely
presented, $\phi:\Sigma\rightarrow Y$ is finite-to-one, and
$D$ is the minimum number of pre-images under $\phi$, then 
\begin{enumerate}
\item $D$ is the minimum number of pre-images under $\pi^0_+$, and
\item there exists an open and dense set $W$ of $Y$ such that if $p\in\mathrm{Per}(Y,f)\cap W$  and $\omega,\omega'\in(\pi^0_+)^{-1}(p)$, then $D=\#\phi^{-1}(p)$ and $\theta(\omega)=\theta(\omega')$.  This follows since $\omega$ and $\omega'$ are canonically associated with points $s$ and $s'$ in $\Sigma$, respectively, where $\phi(s)=\phi(s')$ and we know that such points map to the same point under $\theta$.
\end{enumerate}.
\end{rem}

\begin{cla} If $\Sigma$ is a transitive $1$-step SFT, $(Y,f)$ is finitely
presented, $\phi:\Sigma\rightarrow Y$ is finite-to-one, 
$D$ is the minimum number of pre-images under $\phi$, and  $y\in Y$ is transitive, then $\#(\phi)^{-1}(y)=\#(\pi^0_+)^{-1}(y)=D$ and $\#(\pi_+)^{-1}(y)=1$.
\end{cla}

\noindent{\bf Proof.}
Let $W$ be as in the statement of Claim~\ref{c.a7}.  We know that
$$\min_{y\in Y}\#(\phi)^{-1}(y)=D=\min_{y\in Y}\#(\pi^0_+)^{-1}(y)$$
and for any transitive point $y\in Y$ that $\#(\phi^{-1}(y))=D$ and $\#(\pi^0_+)(y)=D$.  Then the preimages of $y$ for $\pi^0_+$ are the canonically associated points from the preimages of $y$ for $\phi$.  Since all the canonically associated points map to the same point under $\theta$ we know that $\#(\pi_+)^{-1}(y)=1$. $\Box$

\begin{cla}\label{c.transitive}
Suppose $\Sigma$ is a transitive one-step SFT, $Y$ is a transitive Smale spaces, and $\phi:\Sigma\rightarrow Y$ is a finite-to-one factor map.  Let $t,t'\in\Sigma$ such that $\phi(t)=\phi(t')=y$, $t\in W^s(t')$, $t_0\neq t'_0$, and $t_i=t'_i$ for all $i\geq 1$.  If $s$ is a transitive point of $\Sigma$ sufficiently close to $t$ and $\omega\in(\pi^0_+)^{-1}(\phi(s))$, then $t_0\in V^u_{\delta}(\omega)$ and $t'_0\notin V^u_{\delta}(\omega)$.
\end{cla}

\noindent{\bf Proof.}  Since $s$ is transitive we know that words from $\mathcal{W}_{m,n}(w)$ occur infinitely often in the past and in the future for $s$.  From the proof of Claim~\ref{c.magic} we know that 
$$\phi(t)=\phi(t')\in\partial^{t_0}_uR_{t'_0}\cap \partial^{t'_0}_uR_{t_0}.$$  Furthermore, from the proof of Claim~\ref{c.magic} if $s_i=t_i$ for all $-M\leq i\leq M$ and $M$ is sufficiently large, then 
$$\phi([t,s])\notin W^s(\phi(t'), R_{t'_0}).$$
Hence, 
$\phi([t,s])\notin R_{t'_o}$.

If there exists a point $z\in W^u(\phi(s),R_{t_0})\cap R_{t'_0}$, then $$[z,\phi(t)]=[\phi(s),\phi(t)]=\phi[s,t]\in R_{t'_0},$$ a contradiction.  Hence, 
$$W^u(\phi(s),R_{t_0})\cap R_{t'_0}=\emptyset.$$

If $\omega\in(\pi^0_+)^{-1}(\phi(s))$, then $\omega$ is canonically associated with a preimage of $\phi(s)$ under $\phi^{-1}$.  Therefore, $V^u_{\delta}(\omega)$ contains $t_0$ and does not contain $t'_0$.  $\Box$ \\


\noindent{\bf Proof of Lemma~\ref{c.sresolving}} 
We know that $\bar{Y}$ is the unique irreducible component in $\theta(\Sigma_+)$.  This implies that the fiber product of $\Sigma$ and $\bar{Y}$ is transitive.  Hence, $G=G_0$.  Suppose there exists points $(t,y), (t',y)\in G$ where $t\in W^s(t')$.  We may assume that $t_0\neq t'_0$ and $t_i=t'_i$ for all $i\geq 1$.

Let $(s,z)\in G$ be transitive.  By choosing iterates of $(s,z)$ sufficiently close to $(t,y)$ and $(t',y)$ we then know from the characterization of transitive points in Remark~\ref{r.magic} that $t_0,t'_0\in V^u_{\delta}(\omega)$ for all $\omega\in\Sigma_+$ such that $\theta(\omega)=y$.

However, from Claim~\ref{c.transitive} we know that for all iterates of $(s,z)$ sufficiently close to $(t,y)$ that $t'_0\notin V^u_{\delta}(\omega)$ where $\omega$ is contained in the preimage of the iterate of $\tilde{\pi}\pi(s)$ under $\pi^0_+$, a contradiction.  Therefore, $\pi_2$ is $s$-resolving.
$\Box$

\bibliographystyle{plain}
\bibliography{resolve}

\end{document}